\documentclass[12pt]{amsart}
\usepackage[utf8]{inputenc}
\usepackage{style} %moved all style here

\makeatletter
\def\namedlabel#1#2{\begingroup
 #2%
 \def\@currentlabel{#2}%
 \phantomsection\label{#1}\endgroup
}
\makeatother

\title{Generalizations of Knotoids and Spatial Graphs}

\author[C. Adams et al]{Colin Adams}
\address{Department of Mathematics, Williams College,Williamstown, MA 01267}
\email{cadams@williams.edu}
\author[]{Alexandra Bonat}
\address{Department of Mathematics, Williams College,Williamstown, MA 01267} \email{ajb10@williams.edu}
\author[]{Maya Chande}
\address{Department of Mathematics, Princeton University, Princeton, NJ 08544}
\email{mchande@princeton.edu}
\author[]{Joye Chen} 
\address{Department of Mathematics, Princeton University, Princeton, NJ 08544} \email{joyec@princeton.edu}
\author[]{Maxwell Jiang}
\address{Department of Mathematics, M.I.T., Cambridge, MA 02139}
\email{mdjiang@mit.edu}
\author[]{Zachary Romrell}
\address{Department of Mathematics,Williams College, Williamstown, MA 01267} \email{zr3@williams.edu}
\author[]{Daniel Santiago}
\address{Department of Mathematics, M.I.T., Cambridge, MA 02139}
\email{dsantiag@mit.edu}
\author[]{Benjamin Shapiro}
\address{Department of Mathematics,Williams College, Williamstown, MA 01267} \email{bis1@williams.edu}
\author[]{Dora Woodruff}
\address{Department of Mathematics, Harvard University, Cambridge, MA 02138}
\email{dorawoodruff@college.harvard.edu}

\begin{document}

\begin{abstract} In  2010, Turaev introduced knotoids as a variation on knots that replaces the embedding of a circle with the embedding of a closed interval with two endpoints which here we call poles. We define generalized knotoids to allow arbitrarily many poles, intervals, and circles, each pole corresponding to any number of interval endpoints, including zero. This theory subsumes a variety of other related topological objects and introduces some particularly interesting new cases. We explore various analogs of knotoid invariants, including height, index polynomials, bracket polynomials and hyperbolicity. We further generalize to knotoidal graphs,  which are a natural extension of spatial graphs that allow both poles and vertices.
\end{abstract}

\maketitle

%Color coding: Dora - {\color{magenta}magenta}, Daniel - {\color{red}red}, Zach - {\color{violet}violet}, Alex - {\color{orange}orange}, Max - {\color{teal}teal}, Maya - {\color{green}green}, Joye - {\color{purple}purple}, Ben - {\color{pink}pink}

\section{Introduction}\label{introduction section}
Knotoids were introduced by Turaev in \cite{Turaev1} as an extension of classical knot theory. Knotoids are equivalence classes of generic immersions of $[0, 1]$ into $S^2$ considered up to isotopy and Reidemeister moves away from the endpoints. At the endpoints, there is a forbidden fourth move, disallowing an endpoint from moving over/under a strand. Knotoids have been studied recently in various papers, including among many others \cite{winding-poly}, \cite{NIV}, \cite{family-index-poly-knotoid}, \cite{kutluay}. They have also received attention for their potential as protein models, for example in \cite{knotoidproteins} and \cite{knotoidproteins2}.

We introduce \emph{generalized knotoids} and define analogs of knotoid invariants for generalized knotoids. The theory of generalized knotoids in $S^2$ subsumes knot theory and its various extensions, including both planar and spherical (multi-)knotoids as originally defined in \cite{Turaev1}, the multi-linkoids of \cite{multi-linkoid}, long knots as in \cite{Kashaev}, string links as in \cite{HL}, and the polar knots of \cite{polar}. Generalized knotoids are also interesting in their own right, and we include several conjectures concerning their properties.

We also introduce \emph{knotoidal graphs} as a further extension of generalized knotoids. This definition is a natural extension of spatial graphs, bonded knotoids as in \cite{bondedknotoids} and \cite{bondedknotoids2} and graphoids as in \cite{bondedknotoids2} and  \cite{graphoids}. Knotoidal graphs may provide a useful model for proteins that consist of more than one protein molecule bonded together as occurs for hemoglobin and insulin for instance. Finally, we introduce hyperbolicity as an invariant for knotoidal graphs. 

\subsection{Organization}

In Section \ref{generalized knotoids}, we define generalized knotoid theory and describe its relation to various existing extensions of knot theory. The three sections that follow describe analogs of knotoid invariants for generalized knotoids. In Section \ref{height}, we generalize the notion of height and extend a theorem of Kauffman and Gügümcü \cite{gp-min-crossing} regarding height in minimal-crossing diagrams. In Section \ref{index}, we define a family of polynomial invariants for generalized knotoids that extend existing index polynomials for knotoids, and in Section \ref{bracket}, we define a bracket polynomial that generalizes existing bracket polynomial invariants for knotoids.

In Section \ref{knotoidal}, we define knotoidal graphs. We also describe \textit{rail diagrams} for knotoidal graphs, which provide a useful topological point of view. In Section \ref{hyperbolicity}, we use rail diagrams to extend a map defined in \cite{hypknotoids} for spherical knotoids to knotoidal graphs. This map takes knotoidal graphs to spatial graphs in manifolds which are either handlebodies or thickened surfaces. Then we define hyperbolicity and hyperbolic volumes of knotoidal graphs and discuss some applications, particularly to staked knots. 

\subsection*{Acknowledgements} The research was supported by Williams College and NSF Grant DMS-1947438 supporting the SMALL Undergraduate Research Project.

\section{Generalized Knotoids}\label{generalized knotoids}
Let $\Sigma$ denote a closed orientable surface and $G$ a finite graph. We do not require $G$ to be connected or simple, and $G$ may have valency-zero vertices. Let $\tilde G$ denote the disjoint union of $G$ with a finite collection of circles. The edges and circles of $\tilde G$ are called its \emph{constituents}.

A \emph{generalized knotoid diagram} $\diag$ is a generic immersion of $\tilde G$ in $\Sigma$ whose only singularities are transverse double points, called \emph{crossings}, which are labeled with over/undercrossing data. For brevity, we will also use $\diag$ to refer to the image of the immersion. The graph $G$ (resp. $\tilde G)$ is called the \emph{underlying graph} (resp. \emph{underlying looped graph}) of $\diag$. Let $P(\diag)$ denote the set of images of the vertices of $\tilde G$, called the \emph{poles} of $\diag$. The \emph{valency} of a pole $p \in P(\diag)$ is the valency of its corresponding vertex in the underlying graph $G$. A valency zero pole is also called an \emph{isolated pole}.  

Let $E(\diag)$ denote the set of images of the edges of $\tilde G$, called \emph{segment constituents} of $\diag$, and let $L(\diag)$ denote the images of the circles of $\tilde G$, called \emph{loop constituents} of $\diag$. Let $C(\diag) := E(\diag) \cup L(\diag)$ denote the set of \emph{constituents} of $\diag$. A labeling of the vertices of $G$ induces a labeling of the poles of $\diag$. Similarly, a labeling (resp. an orientation) of the constituents of $\tilde G$ induces a labeling (resp. orientation) of the constituents of $\diag$.

\begin{figure}[htbp]
\begin{center}
\includegraphics[width=9cm]{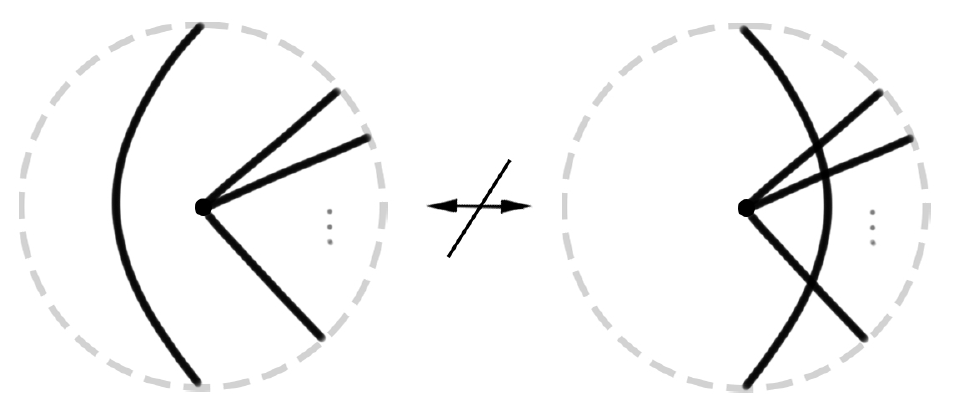}

\includegraphics[width=9cm]{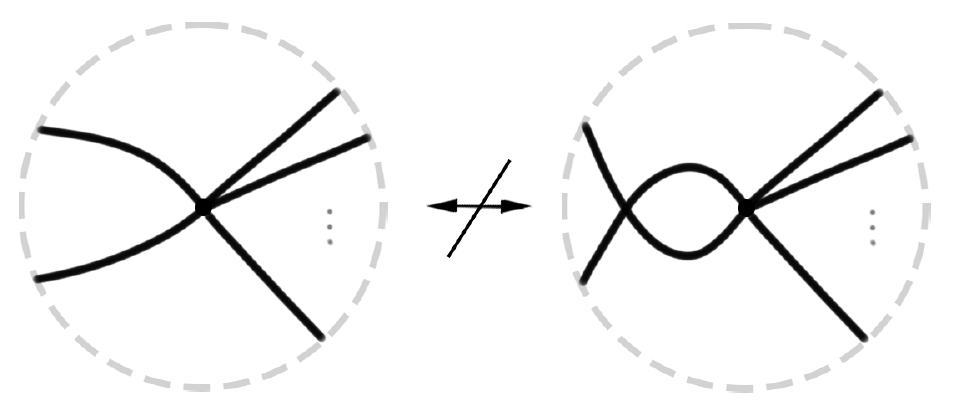}
\caption{The pole slide move (top) and the pole twist move (bottom) are forbidden, irrespective of the crossing data chosen for the diagrams on the right.}
\label{forbidden}
\end{center}
\end{figure}

We introduce an equivalence relation on generalized knotoid diagrams in $\Sigma$ generated by ambient isotopy and the three standard Reidemeiester moves away from the poles. Observe that it is forbidden for a constituent to pass through a pole and for ``twists'' to be created or destroyed near a pole. See Figure \ref{forbidden}. We remark that the forbidden twist move of Figure \ref{forbidden} parallels the forbidden move for vertices in a rigid-vertex spatial graph.

A \emph{generalized knotoid} $\kappa$ is an equivalence class of generalized knotoid diagrams. It is clear that equivalency respects the number of poles and constituents as well as their associated data; in particular, we may speak of labeled and/or oriented generalized knotoids as equivalence classes of labeled and/or oriented generalized knotoid diagrams. Note that allowing twisting near poles yields a different theory of generalized knotoids, which we do not consider in this paper. However, some of the invariants we define later are able to distinguish between generalized knotoids which differ only by twists near poles and are otherwise identical. 

An equivalent, topological viewpoint is to define a generalized knotoid diagram as a generic immersion of $\Tilde{G}$, where all vertices of $G$ have valency at least one, in a compact orientable surface $\Sigma$ with or without boundary. Given such a diagram on a compact surface, we may recover a diagram on a closed surface by capping off each boundary component with a disk, then collapsing each disk to a point representing an isolated pole. Conversely, given a generalized knotoid diagram on a closed surface with isolated poles, we may remove an open disk neighborhood of each isolated pole. That is, the theory of generalized knotoids on closed surfaces with isolated poles is equivalent to the theory of generalized knotoids on compact surfaces without isolated poles. We rely on both viewpoints, especially in Section \ref{hyperbolicity} where the topological viewpoint proves useful. 

We denote the case $\Sigma = S^2$ by \emph{classical generalized knotoid theory}. In the remainder of this paper, a generalized knotoid is assumed to be classical unless otherwise stated.

\begin{figure}
\captionsetup[subfigure]{labelfont=rm}
    \centering
    \begin{subfigure}{0.3\textwidth}
    \centering
    \includegraphics[width=3cm]{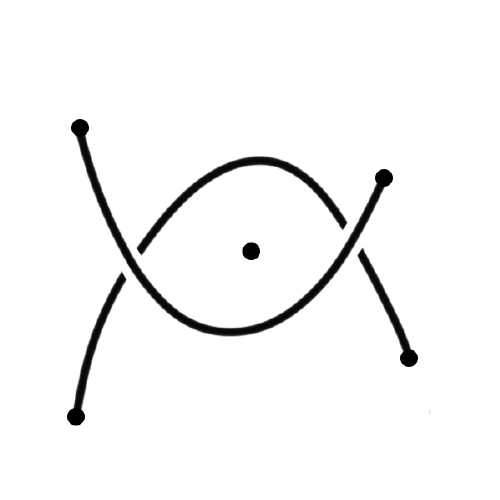}
    \caption{}
    \end{subfigure}
    \begin{subfigure}{0.3\textwidth}
    \centering
    \includegraphics[width=3cm]{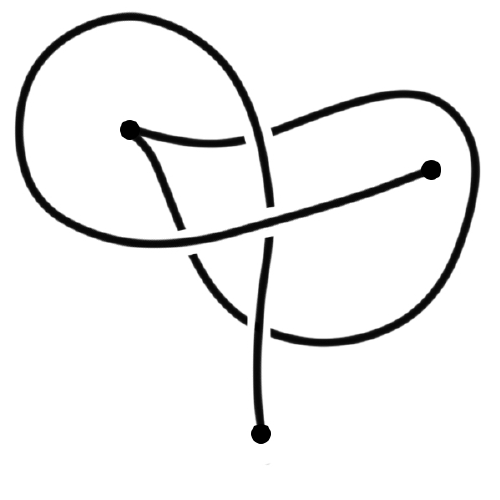}
    \caption{}
    \label{bee}
    \end{subfigure}
    \begin{subfigure}{0.3\textwidth}
    \centering
    \includegraphics[width=3cm]{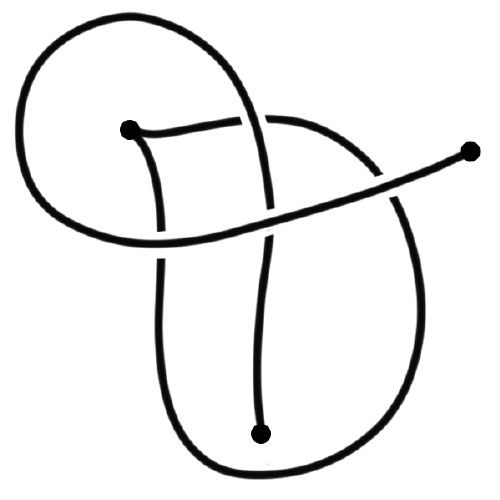}
    \caption{}
    \label{cee}
    \end{subfigure}
    \begin{subfigure}{0.3\textwidth}
    \centering
    \includegraphics[width=3cm]{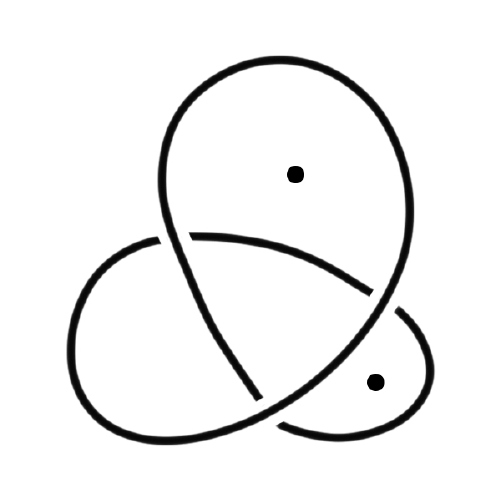}
    \caption{}
    \end{subfigure}
    \begin{subfigure}{0.3\textwidth}
    \centering
    \includegraphics[width=3cm]{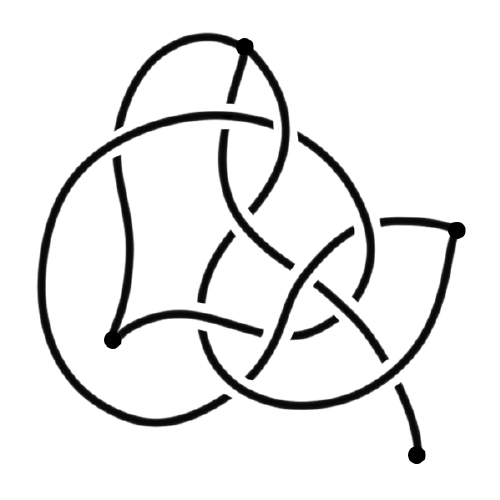}
    \caption{}
    \end{subfigure}
    \begin{subfigure}{0.3\textwidth}
    \centering
    \includegraphics[width=3cm]{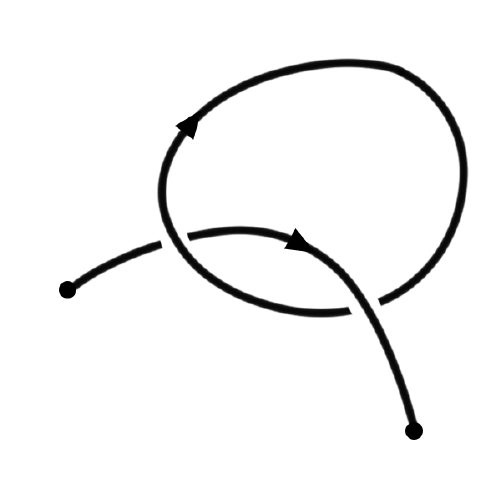}
    \caption{}
    \label{eff}
    \end{subfigure}
    \caption{Generalized knotoid diagrams.}
    \label{genknotexamples}
\end{figure}

\begin{example}
See Figure \ref{genknotexamples} for some diagrams of generalized knotoids. The reader may verify that Figures \ref{bee} and \ref{cee} represent equivalent generalized knotoids. Figure \ref{eff} represents an oriented generalized knotoid.
\end{example}
\begin{example}\label{spherical_example}
A spherical knotoid is a generalized knotoid with two poles and a single segment constituent between them. 

\end{example}

\begin{example}\label{planar_example}
A planar knotoid is a generalized knotoid with three poles and a single segment constituent connecting two of them. The equivalency of the two notions is obtained by identifying the isolated pole of the generalized knotoid with the point at infinity of the planar knotoid diagram. This example provides a useful viewpoint for a number of invariants defined for planar knotoids, and it motivates parts of the constructions of invariants for generalized knotoids in Sections \ref{index} and \ref{bracket}.
\end{example}

\begin{example}
A classical knot or link is a generalized knotoid with no poles.
\end{example}

\begin{example}
The multi-linkoids of \cite{multi-linkoid} are generalized knotoids whose poles all have valency one.
\end{example}

\begin{example}
For a positive integer $n$, an \emph{$n$-polar knot} \cite{polar} is a generalized knotoid with no loop constituents and whose underlying graph is the cycle graph $C_n$. An $n$-polar knot diagram looks like a classical knot diagram with $n$ poles placed on the knot, away from the crossings. Long knots \cite{Kashaev} can be thought of as $1$-polar knots, and the equivalency of the two notions is obtained by identifying the single pole of a polar knot diagram with the point at infinity in a long knot diagram.
\end{example}

\begin{figure}[htbp]
    \captionsetup[subfigure]{labelfont=rm}
    \begin{subfigure}{0.4\textwidth}
    \centering
    \includegraphics[]{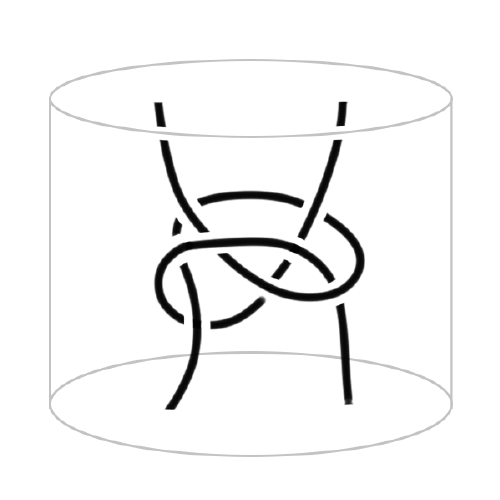}
    \caption{}
    \label{stringlink}
    \end{subfigure}
    \begin{subfigure}{0.4\textwidth}
    \centering
    \includegraphics[]{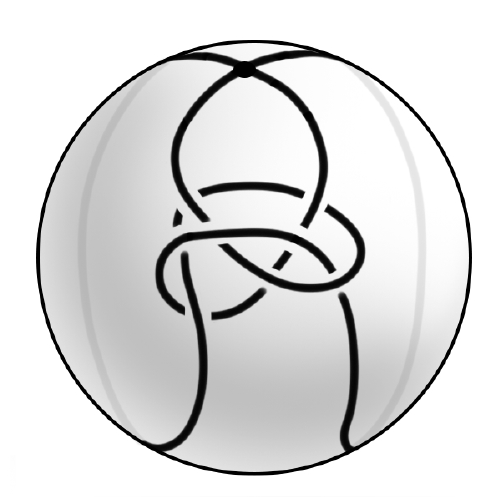}
    \caption{}
    \label{stringlinksphere}
    \end{subfigure}
    \caption{The string link in Figure \ref{stringlink} corresponds to the generalized knotoid in Figure \ref{stringlinksphere}. The fact that Reidemeister moves must occur away from the pole captures the fact that the endpoints of the string link's segments are fixed.}
    \label{stringlinkexample}
\end{figure}

\begin{example} A \emph{string link} with $n$ strings is defined in \cite{HL} as the embedding of a finite set of closed intervals $I_1, I_2, \dots, I_n$ in $D \times I$ where $D$ is the unit disk in the $xy$-plane, such that the initial point of $I_i$ is sent to $(x_i, 0,0)$,  and the final point is sent to  $(x_i, 0,1)$,  where $-1 < x_1 < x_2 < \dots < x_n < 1$. This is defined up to ambient isotopy of the strings with fixed endpoints.

We represent a string link as a generalized knotoid with a single pole of valency $2n$. See Figure \ref{stringlinkexample}. If we thicken the sphere, we may consider the pole as the removal of a neighborhood of a single vertical line in $S^2 \times I$ (bringing the ambient space to $D\times I$) and the constituents as embedded segments with endpoints on the circle boundary of $D \times \{1/2\}$. (See Section \ref{railsection} for a discussion on \emph{rail diagrams}.) The legs of the strings appear in clockwise order $\ell_1, \ell_2, \dots, \ell_n$ around  the circle boundary of $D \times \{1/2\}$ followed by the heads of the strings in clockwise order $h_n, h_{n-1}, \dots, h_1$. Since the pole twist move is disallowed, these endpoints are fixed in their order around the circle, and we recover the notion of a string link. 
\end{example}

\begin{example} \label{staked-links}
We define \emph{staked links} on a surface $\Sigma$ to be the class of generalized knotoids on $\Sigma$ with no segment constituents and at least one pole. A diagram of a staked link looks like a link diagram on $\Sigma$ with isolated poles placed in the regions of the link's complement in $\Sigma$. An example appears in Figure \ref{genknotexamples}(d).
\end{example}

\section{Height for Generalized Knotoids}\label{height}
We review the notion of height for knotoids, defined originally as \emph{complexity} in \cite{Turaev1}. Let $\diag$ be a knotoid diagram on $\Sigma$, and let $\alpha$ be an embedded arc in $\Sigma$ connecting the endpoints that intersects $\diag$ transversely and away from the crossings. We call $\alpha$ a \emph{shortcut}. The \emph{diagram height} $h(\diag)$ is the number of intersections of $\alpha$ with $\diag$ not including the endpoints, minimized over all shortcuts $\alpha$. The \emph{height} of a knotoid $k$ is the minimum of $h(D)$ over all diagrams $\diag$ representing $k$, and a diagram $\diag$ attaining this minimum \emph{realizes the height}.

This definition extends to generalized knotoids in a natural way. Let $\diag$ be a pole-labeled generalized knotoid diagram on $\Sigma$ and let $A,B \in P(\diag)$ be distinct poles. A \emph{shortcut} $\alpha$ between $A$ and $B$ is an embedded arc in $\Sigma$ with endpoints at $A$ and $B$ that intersects $\diag$ transversely and away from crossings and poles. (If $A=B$, we define a shortcut in the same way with ``embedded arc'' replaced by ``simple closed curve''.) The \emph{diagram height} between $A$ and $B$, denoted $h_{\diag}(A,B)$, is the number of intersections of $\alpha$ with $\diag$ not including $A$ and $B$, minimized over all shortcuts $\alpha$ between $A$ and $B$. An $\alpha$ attaining this minimum is a \emph{minimal shortcut}. Given a generalized knotoid $\kappa$ and poles $A,B \in P(\kappa)$, the \emph{height} between $A$ and $B$, denoted $h_\kappa(A,B)$, is the minimum of $h_{\diag}(A,B)$ over all diagrams $\diag$ representing $K$. A diagram $\mathcal D$ attaining this minimum \emph{realizes the height between $A$ and $B$}. Given an ordering on the poles $P_1, \dots, P_n$ of $\kappa$, the \emph{height spectrum} of $\kappa$ is the symmetric matrix $(h(P_i, P_j))_{ij}$.

Theorem \ref{originalHeightTheorem} was conjectured by Turaev in \cite{Turaev1} and proved by Gügümcü and Kauffman in \cite{gp-min-crossing}.
\begin{theorem}[\cite{gp-min-crossing}, Theorem 13]\label{originalHeightTheorem}
Let $k$ be a spherical knotoid with $h(k) = 0$. Then any minimal-crossing diagram $\mathcal D$ of $k$ realizes the height.
\end{theorem}
The proof of Theorem \ref{originalHeightTheorem} given in \cite{gp-min-crossing} considers the \emph{virtual closure} of a knotoid and applies results from virtual knot theory derived from the method of \emph{parity projection} \cite{manturovparity}. In Theorem \ref{heightTheorem}, we extend Theorem \ref{originalHeightTheorem} to generalized knotoids. The proof of Theorem \ref{heightTheorem} adapts a topological approach to parity via colorings described in \cite{boden-rushworth}.

\begin{theorem}\label{heightTheorem}
Let $\kappa$ be a classical generalized knotoid and suppose poles $A,B \in P(\kappa)$ satisfy $h_\kappa(A,B) = 0$. Then any minimal-crossing diagram $\diag$ of $\kappa$ realizes the height between $A$ and $B$.
\end{theorem}

\begin{proof}
    If $A = B$, the result is trivial, so assume $A \neq B$. Fix a diagram $\diag_0$ such that $h_{\diag_0}(A,B) = 0$, and let $\alpha$ be a minimal shortcut between $A$ and $B$. Suppose for contradiction's sake that there exists a minimal-crossing diagram $\diag$ with $h_{\diag}(A,B) > 0$. We may assume that there is a sequence of diagrams $\diag_0 \to \diag_1 \to \dots \to \diag_n = \diag$ where for each $0 \le i < n$, the diagrams $\diag_i$ and $\diag_{i+1}$ are related by a Reidemeister move. Recall that Reidemeister moves take place away from poles, so we may assume that each pole $Q$ of $\diag_0$ has a neighborhood $N_Q$ that is unchanged as the moves occur. We declare that the shortcut $\alpha$ is fixed across all the diagrams $\diag_i$, and we may assume that $\alpha$ intersects each $\diag_i$ transversely and away from crossings.
    \begin{claim}\label{even-intersection-Claim}
        The number of intersections of $\alpha$ with a constituent of $\diag_i$, excluding potential intersections at poles, is even.
    \end{claim}
    \begin{proof}
        Let $e \in C(\kappa)$ be a constituent. For each $0 \le j < n$, let $\#(\alpha \cap e)_j$ denote the number of intersections of $\alpha$ with $e$ in the diagram $\diag_j$, excluding poles. We show that $\#(\alpha \cap e)_j$ and $\#(\alpha \cap e)_{j+1}$ have the same parity, and the conclusion follows from the assumption that $\#(\alpha \cap e)_0 = 0$.
        
        Let $D$ be the disk region where the Reidemeister move $\diag_j \to \diag_{j+1}$ is applied. For each $t \in \{1, 2, 3\}$, a type $t$ Reidemeister move has $t$ \emph{participating strands} $s_1, \dots s_t$, some of which may belong to $e$. Each $s_\ell$ meets $\partial D$ at two \emph{endpoints} that are fixed by the Reidemeister move.
        
        If $\alpha$ does not intersect $D$, then the conclusion is clear. Otherwise, we may assume that $\alpha$ meets $\partial D$ at a finite number of points $P_1, P_2, \dots, P_k$ and $Q_1, Q_2, \dots, Q_k$, distinct from the endpoints of the participating strands, such that each pair $P_m$ and $Q_m$ is connected by a non-self-intersecting arc $\alpha_m \subset \alpha$ contained in the interior of $D$. Each $\alpha_m$ divides $D$ into two regions, and $\#(\alpha_m \cap s_\ell)_j$ is odd if the endpoints of $s_\ell$ are in opposite regions of $D$ and even otherwise. Since the Reidemeister move fixes the endpoints, it follows that $\#(\alpha_m \cap s_\ell)_j$ and $\#(\alpha_m \cap s_\ell)_{j+1}$ have the same parity for each $m$ and $\ell$, and the conclusion follows. See Figure \ref{RTwoShortcut}.
    \end{proof}
    \begin{figure}[htbp]
        \centering
        \includegraphics[width=8cm]{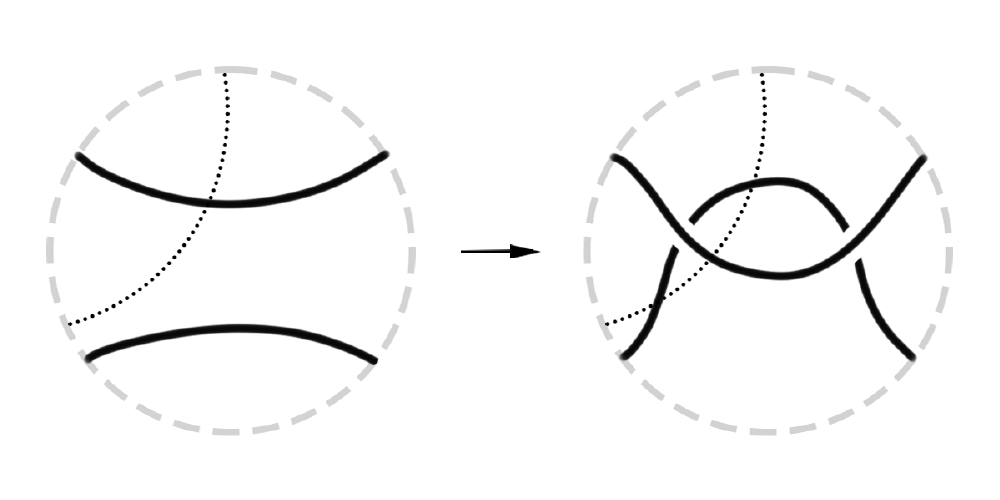}
        \caption{A Reidemeister move occuring in the disk region $D$, whose boundary is indicated by dashed line. The shortcut $\alpha$ is indicated with dotted line. The parity of the number of intersections of $\alpha$ with each participating strand is preserved.}
        \label{RTwoShortcut}
    \end{figure}

    Given a constituent $e$ in a diagram $\diag_i$, an \emph{$\alpha$-coloring} of $e$ is an assignment of black or gray to each point of $e$ such that the color of $e$ changes precisely at the points where it intersects $\alpha$. It follows from Claim \ref{even-intersection-Claim} that an $\alpha$-coloring exists for each $e$. We define an \emph{$\alpha$-coloring} of the diagram $\diag_i$ to be a choice of $\alpha$-coloring for each of its constituents such that each segment constituent is black in a neighborhood of its endpoint poles. We inductively fix an $\alpha$-coloring $\mathcal C_i$ for each $\diag_i$ as follows. We declare $\mathcal C_0$ to be the $\alpha$-coloring in which all constituents of $\diag_0$ are colored entirely black. Once $\mathcal C_i$ is chosen, the Reidemeister move $\diag_i \to \diag_{i+1}$ uniquely determines $\mathcal C_{i+1}$ if we require $\mathcal C_i$ and $\mathcal C_{i+1}$ to agree outside the disk region where the Reidemeister move is applied. See Figure \ref{alphaColoring}. Under the $\alpha$-coloring $\mathcal C_i$, a crossing $c$ in $\diag_i$ is assigned a color from each of the constituents meeting at $c$. We say $c$ is \emph{even} if these colors agree and \emph{odd} if they differ.
    
    \begin{figure}[htbp]
    \captionsetup[subfigure]{labelfont=rm}
        
        \begin{subfigure}{0.3\textwidth}
        \centering
        \includegraphics[width=3cm]{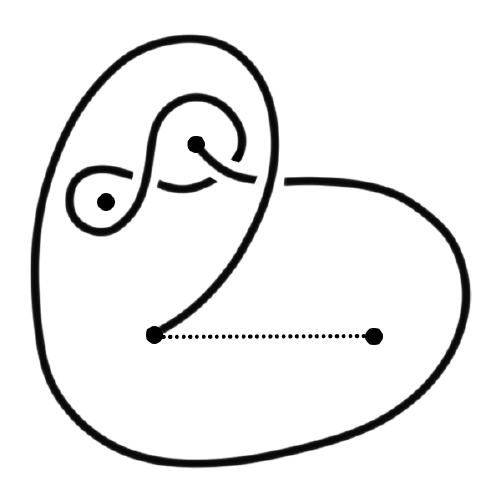}
        \caption{}
        \label{diagZero}
        \end{subfigure}
        \begin{subfigure}{0.3\textwidth}
        \centering
        \includegraphics[width=3cm]{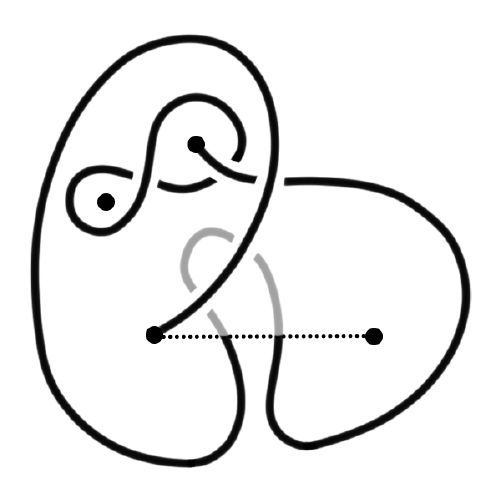}
        \caption{}
        \label{diagOne}
        \end{subfigure}
        \begin{subfigure}{0.3\textwidth}
        \centering
        \includegraphics[width=3cm]{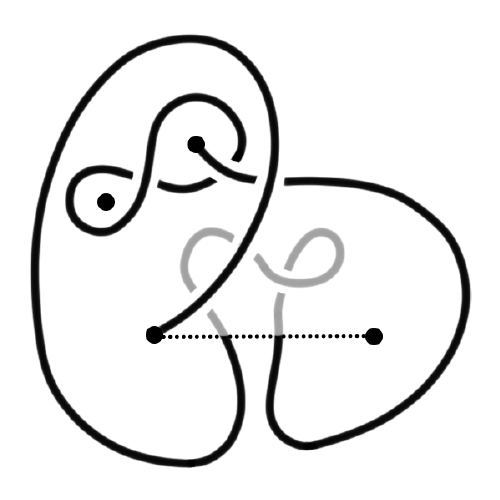}
        \caption{}
        \label{diagTwo}
        \end{subfigure}
        
        \caption{Figures \ref{diagZero}, \ref{diagOne}, \ref{diagTwo} show an example sequence of diagrams $\diag_0$, $\diag_1$, $\diag_2$ equipped with their $\alpha$-colorings. The shortcut $\alpha$ is indicated with dotted line.}
        \label{alphaColoring}
    \end{figure}
    
    Our next aim is to show that the diagram $\diag = \diag_n$ has at least one odd crossing. To do this, we introduce some notation.
    
    Define a \emph{region} to be a connected component of the complement $S^2 \setminus \diag$. Observe that a region need not be an open disk, and its boundary need not be connected. For a given region $R$, define $\partial_0 R$ to be the subset of $\partial R \subset \diag$ consisting of points $x$ such that every neighborhood of $x$ intersects at least two distinct regions. In particular, the shortcut $\alpha$ intersects $\partial_0 R$ each time it enters or exits a region $R$. Conversely, each intersection of $\alpha$ with $\partial_0 R$ corresponds to an entering or exiting of $R$.
    
    Define a \emph{looped graph} to be the disjoint union of a finite graph with a collection of circles. A looped graph is \emph{Eulerian} if the edges of its graph portion can be partitioned into cycles (i.e. closed paths). It is a classical result that a looped graph is Eulerian if and only if all of its vertices have even valency.
    
    We use the term \emph{singularity} to refer to a pole or crossing of $\diag$. We use the term \emph{border} to refer to a portion of a constituent $e$ that connects consecutive singularities along $e$, or all of $e$ if $e$ is a loop constituent with no singularities. Every singularity $s$ has a disk neighborhood $N_s$ that contains no other singularities. In particular, $\diag$ meets $N_s$ at $m$ radii of $N_s$, where $m$ is the number of borders meeting at $s$.
    
    \begin{claim}\label{eulerian}
        The set $\partial_0R$ forms an Eulerian looped graph whose vertices are singularities and whose edges and circles are borders.
    \end{claim}
    \begin{proof}
        First note that if any point of a border is in $\partial_0R$, then the entire border (including singularities at its endpoints) is contained in $\partial_0R$. It now suffices to show that each singularity $s \in \partial_0R$ has even valency. Consider $N_s$ and $m$ as defined above. The diagram $\diag$ divides $N_s$ into $m$ connected components belonging to (not necessarily distinct) regions $R_1, \dots, R_m$, some of which are equal to $R$. The radii belonging to $\partial_0R$ are those that separate an $R$-region and a non-$R$-region. There are an even number of such radii, and the conclusion follows.
    \end{proof}
    \begin{claim}\label{oddCrossings}
        The diagram $\diag = \diag_n$ has at least one odd crossing.
    \end{claim}
    \begin{proof}
        For sufficiently small choices of neighborhoods $N_A$ and $N_B$, there is a unique \emph{starting region} $R_A$ satisfying $N_A \cap \alpha \setminus \{A\} \subset R_A$. See Figure \ref{startingRegion}. Similarly, there is a unique \emph{ending region} $R_B$ satisfying $N_B \cap \alpha \setminus \{B\} \subset R_B$.  By Claim \ref{eulerian}, we may view $\partial_0R_A$ an Eulerian looped graph and partition its edges into cycles. By the assumption that $h_\diag(A,B) \neq 0$, we have $R_A \neq R_B$. Thus, $\alpha$ intersects $\partial_0 R_A$ an odd number of times, and it follows that $\alpha$ intersects some cycle or circle $C \subset \partial_0R_A$ an odd number times.
        
        Consider the $\alpha$-coloring of $\diag$ restricted to $C$. Recall that constituents must be black near poles, so $C$ does not change color at poles. Thus, each point of $C$ at which its color changes falls into one of two categories:
        \begin{enumerate}[(i)]
            \item An intersection of $C \setminus \{A, B\}$ with $\alpha$. (Every such intersection yields a color change.)
            \item A point of $C$ at which $\diag$ has an odd crossing.
        \end{enumerate}
        As we make a full traversal around $C$, the color changes an even number of times. An odd number of color changes fall under (i), so an odd number of color changes fall under (ii). In particular, $\diag$ has at least one odd crossing.
    \end{proof}
    
    \begin{figure}[htbp]
    \captionsetup[subfigure]{labelfont=rm}
        \centering
        \begin{subfigure}{0.4\textwidth}
        \centering
        \includegraphics[width=4cm]{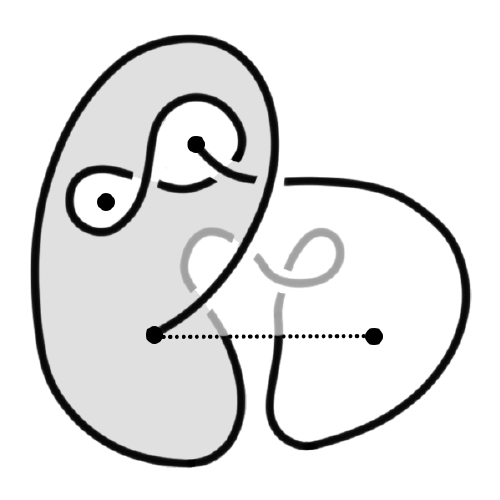}
        \caption{}
        \label{startingRegion}
        \end{subfigure}
        \begin{subfigure}{0.4\textwidth}
        \centering
        \includegraphics[width=4cm]{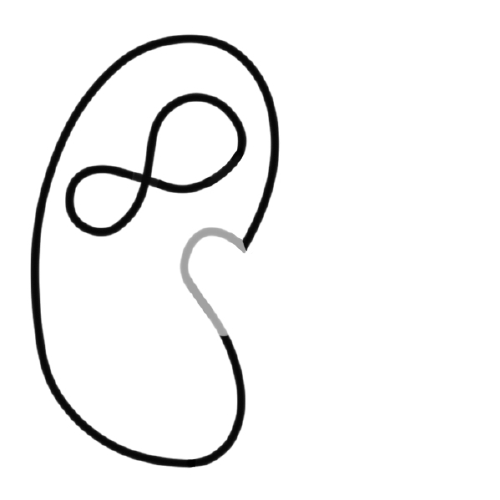}
        \caption{}
        \label{boundaryZero}
        \end{subfigure}
        \caption{The diagram from Figure \ref{diagTwo}. Assuming $A$ is the left endpoint pole of the shortcut $\alpha$, Figure \ref{startingRegion} shows the starting region $R_A$, shaded. Figure \ref{boundaryZero} shows $\partial_0R_A$. Note that one of the color changes on the outer cycle of $\partial_0R_A$ occurs at an odd crossing.}
    \end{figure}
    
    Fix a double branched cover $p \colon S^2 \to S^2$ branched over $\alpha$. Let $e \in C(\diag_i)$ be a constituent. (If $e$ has an endpoint at $A$ or $B$, replace $e$ by $e \setminus \{A, B\}$.) The preimage $p^{-1}(e)$ has two connected components $e', e''$, each homeomorphic to $e$. The preimage $p^{-1}(\diag_i)$ defines a generalized knotoid diagram on the double cover sphere with poles at the points of $p^{-1}(P(\diag_i))$ and constituents $e', e''$ for each $e \in C(\diag_i)$. Its over/under crossing data is inherited from $\diag_i$ in the natural way.
    
    The preimage $p^{-1}(\alpha)$ is a circle, and $S^2 \setminus p^{-1}(\alpha)$ is the disjoint union of two open disks $H_1$ and $H_2$. Let $p_1 \colon H_1 \to S^2 \setminus \alpha$ and $p_2 \colon H_2 \to S^2 \setminus \alpha$ be the restrictions of $p$ to $H_1$ and $H_2$. Note that $p_1$ and $p_2$ induce bijections $p_1 \colon p^{-1}(e) \cap H_1 \to e$ and $p_2 \colon p^{-1}(e) \cap H_2 \to e$, and recall that $e$ has an $\alpha$-coloring dictated by $\mathcal C_i$. We color each point of $p^{-1}(e)$ such that the bijections induced by $p_1$ and $p_2$ are color-preserving and color-reversing, respectively. Doing so for each $e \in C(\diag_i)$ yields the \emph{induced $\alpha$-coloring} of $p^{-1}(\diag_i)$. See Figure \ref{doubleCover}. Using the induced $\alpha$-colorings of the preimages $p^{-1}(e)$, we define \emph{even} (resp. \emph{odd}) crossings of $p^{-1}(\diag_i)$ as crossings where the colors agree (resp. differ). Observe that $p$ defines a parity-preserving, $2$-to-$1$ map from the crossings of $p^{-1}(\diag_i)$ to the crossings of $\diag_i$.
    
    \begin{figure}[htbp]
        \centering
        \includegraphics[width=12cm]{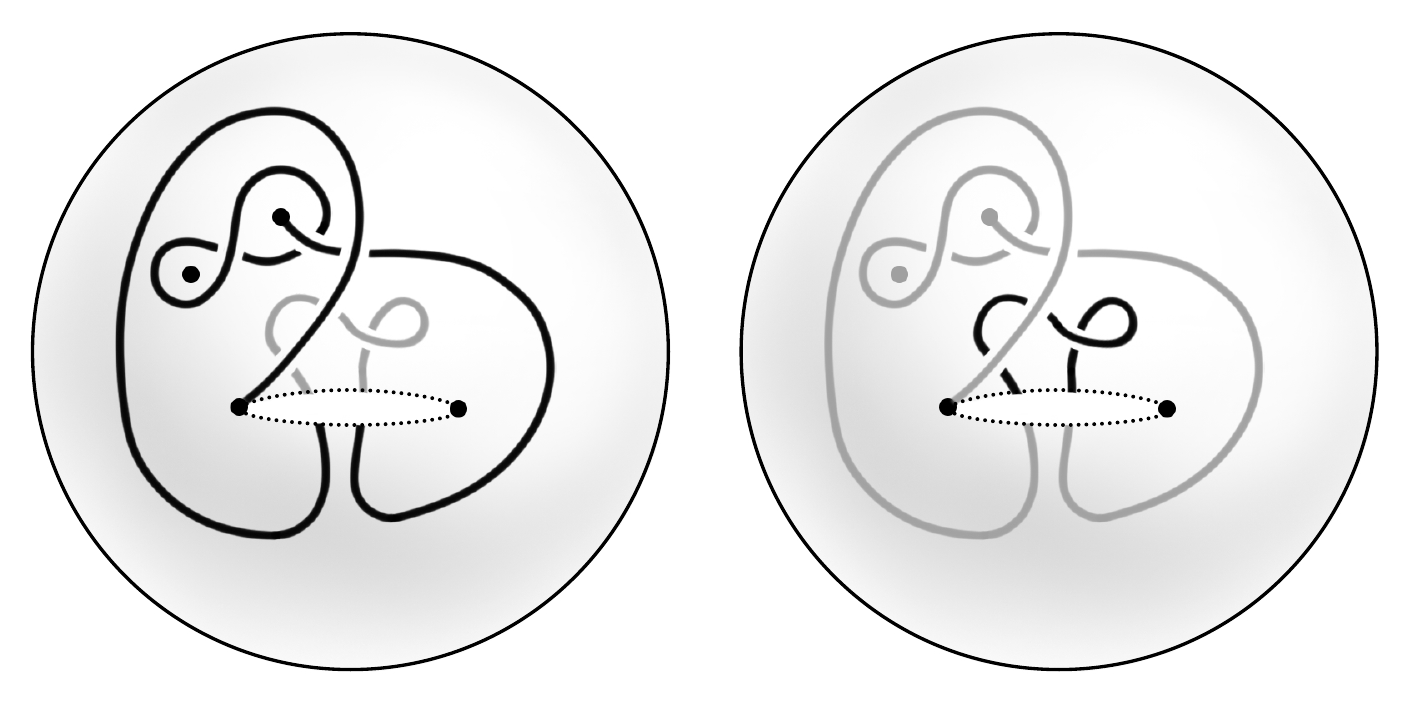}
        \caption{One possible visualization of the branched double cover. On the left is $H_1$, on the right is $H_2$, and they are glued together along their boundaries (the dotted circles obtained by cutting a sphere along $\alpha$) such that the top and bottom arcs on the left are identified with the bottom and top arcs on the right. The preimage $p^{-1}(\diag)$ and its induced $\alpha$-coloring are shown for the diagram $\diag$ from Figure \ref{diagTwo}.}
        \label{doubleCover}
    \end{figure}
    
    \begin{claim}\label{monochrome-lift}
        Under the induced $\alpha$-coloring of $p^{-1}(\diag_i)$, one of $e'$ and $e''$ is entirely black and the other is entirely gray. Moreover, if $e$ has an endpoint at a pole $Q \not \in \{A,B\}$, then the black component of $p^{-1}(e)$ has a corresponding endpoint at $p_1^{-1}(Q) \in H_1$.
    \end{claim}
    \begin{proof}
        It is clear from the definition of the induced $\alpha$-coloring that $p^{-1}(e)$ has both black and gray points. Each time $e$ crosses $\alpha$ and changes color, the lift $e_1$ crosses from $H_1$ to $H_2$ or vice versa. Since $p_1$ is color-preserving and $p_2$ is color-reversing, $e_1$ does not change color. It follows that $e_1$ is monochromatic. Similarly, $e_2$ is monochromatic. This proves the first statement.
        
        Now suppose $e$ has an endpoint at a pole $Q \not \in \{A, B\}$. By definition of $\alpha$-coloring, there is a neighborhood $N_Q$ of $Q$ such that $\mathcal C_i$ assigns black to $e \cap N_Q$. Since $p_1$ is color-preserving, there is a neighborhood of $p_1^{-1}(Q)$ in which $p^{-1}(e)$ is black. This proves the second statement.
    \end{proof}
    In light of Claim \ref{monochrome-lift}, we let $e'$ denote the component of $p^{-1}(e)$ colored black and define a generalized knotoid diagram $\diag_i' \subset p^{-1}(\diag_i)$ as follows. The diagram $\diag_i'$ has poles at $p^{-1}(A)$, $p^{-1}(B)$, and $p_1^{-1}(Q) \in H_1$ for each $Q \in P(\diag_i) \setminus \{A,B\}$, and it has the constituent $e'$ for each $e \in C(\diag_i)$. Its over/under crossing data is inherited from $\diag_i$ in the natural way. The black portions of Figure \ref{doubleCover} show $\diag'$, where $\diag = \diag_2$ is the diagram from Figure \ref{diagTwo}.
    \begin{claim}\label{evenCrossings}
        The crossings of $\diag_i'$ are in bijection with the even crossings of $\diag_i$.
    \end{claim}
    \begin{proof}
        Recall that, as a map from the crossings of $p^{-1}(\diag_i)$ to the crossings of $\diag_i$, the map $p$ preserves parity. Each crossing $c'$ of $\diag_i'$ is even (black-black), so $p(c')$ is an even crossing of $\diag_i$. Thus $p$ defines map from the crossings of $\diag_i'$ to the even crossings of $\diag_i$. For any even crossing $c$ in $\diag_i$, the preimage $p^{-1}(c) = \{p_1^{-1}(c), p_2^{-1}(c)\}$ consists of two crossings, one black-black and one gray-gray. Only the black-black crossing is a crossing of $\diag_i'$. It follows that $p$ gives the desired bijection.
    \end{proof}
    
    Let $\sim$ denote the equivalence generated by ambient isotopies and Reidemeister moves.
    
    \begin{claim}\label{reidemeisterLift}
        For each $0 \le i < n$, we have $\diag_i' \sim \diag_{i+1}'$.
    \end{claim}
    \begin{proof}
        Let $D$ be the disk region where the type $t$ Reidemeister move $\diag_i \to \diag_{i+1}$ is applied and let $s_1, \dots, s_t$ be the participating strands, with endpoints on $\partial D$, as in the proof of Claim \ref{even-intersection-Claim}. The preimage $p^{-1}(D)$ is the disjoint union of two disk regions $D_1$ and $D_2$, and the move $\diag_i \to \diag_{i+1}$ lifts to a move $p^{-1}(\diag_i) \to p^{-1}(\diag_{i+1})$ given by a type $t$ Reidemeister move in each of $D_1$ and $D_2$. By Claim \ref{monochrome-lift}, each participating strand $s_\ell$ lifts to a monochromatic strand $s_{\ell, 1} \subset D_1$ and a monochromatic strand $s_{\ell, 2} \subset D_2$ (of the opposite color). Since $\diag_i \to \diag_{i+1}$ preserves the colors of the endpoints, the move $p^{-1}(\diag_i) \to p^{-1}(\diag_{i+1})$ also preserves the colors of the endpoints and thus preserves the colors of $s_{\ell, 1}$ and $s_{\ell, 2}$. It follows that restricting the move $p^{-1}(\diag_i) \to p^{-1}(\diag_{i+1})$ to only the black participating strands yields either a Reidemeister move or an ambient isotopy that takes $\diag_i'$ to $\diag_{i+1}'$. See Figure \ref{liftRThree}.
    \end{proof}
    
    \begin{figure}[htbp]
    \captionsetup[subfigure]{labelfont=rm}
    \centering
    \begin{subfigure}{0.45\textwidth}
    \centering
    \includegraphics[width=5.7cm]{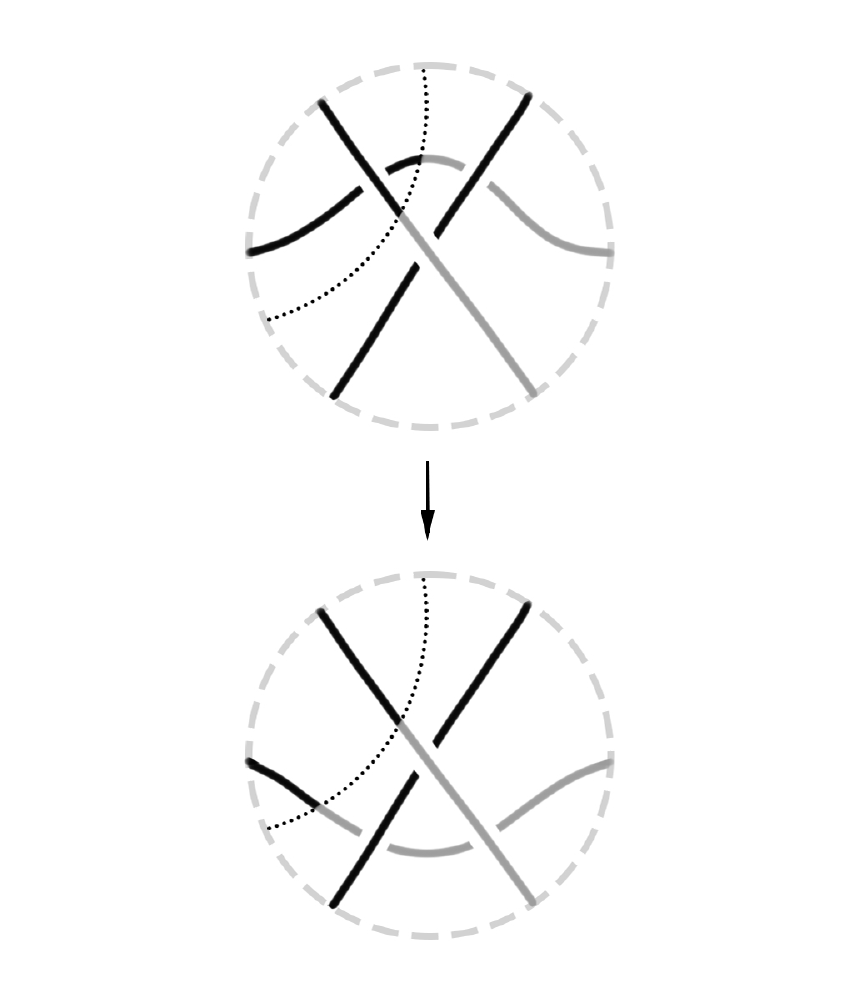}
    \caption{}
    \label{coloredRThree}
    \end{subfigure}
    \begin{subfigure}{0.45\textwidth}
    \centering
    \includegraphics[width=5.7cm]{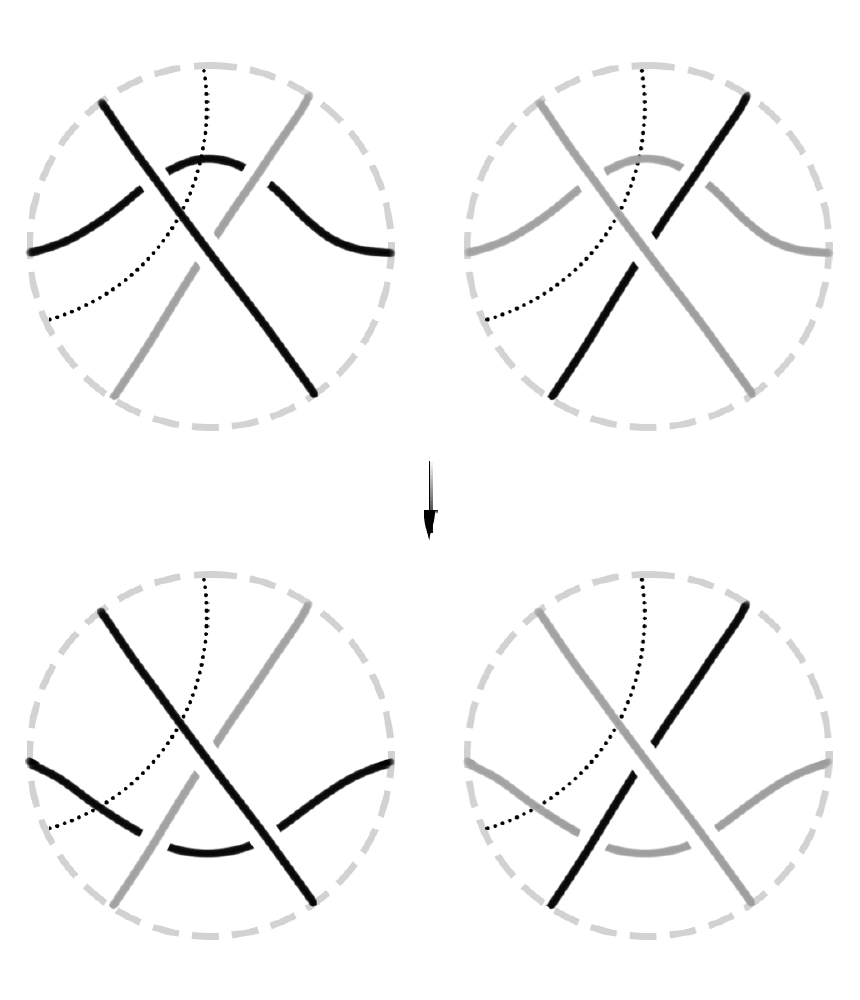}
    \caption{}
    \label{liftedRThree}
    \end{subfigure}
    \caption{A Reidemeister move $\diag_i \to \diag_{i+1}$ (Figure \ref{coloredRThree}) lifts to a move $p^{-1}(\diag_i) \to p^{-1}(\diag_{i+1})$ in two disk regions (Figure \ref{liftedRThree}). Restricting to only the black strands in Figure \ref{liftedRThree} yields a move $\diag_i' \to \diag_{i+1}'$. In this example it is an ambient isotopy.}
    \label{liftRThree}
\end{figure}
    We finish the proof of Theorem \ref{heightTheorem}. Since $D_0$ does not intersect $\alpha$ except at $A$ and $B$, it follows by construction that $\diag_0'$ is contained in $H_1 \cup \{p^{-1}(A), p^{-1}(B)\}$, and that $\diag_0$ and $\diag_0'$ are ambiently isotopic diagrams of $\kappa$. See Figure \ref{doublecoverdiagzero}. We now have
    \[\diag_n \sim \diag_0 \sim \diag_0' \sim \diag_n',\]
    where the last equivalence follows from Claim \ref{reidemeisterLift}. Hence $\diag_n$ and $\diag_n'$ are both diagrams of $\kappa$. However, Claim \ref{oddCrossings} and Claim \ref{evenCrossings} together imply that $\diag_n'$ has fewer crossings than $\diag_n$. This contradicts the assumption that $\diag_n$ is a minimal-crossing diagram of $\kappa$.
    \begin{figure}[htbp]
        \centering
        \includegraphics[width=12cm]{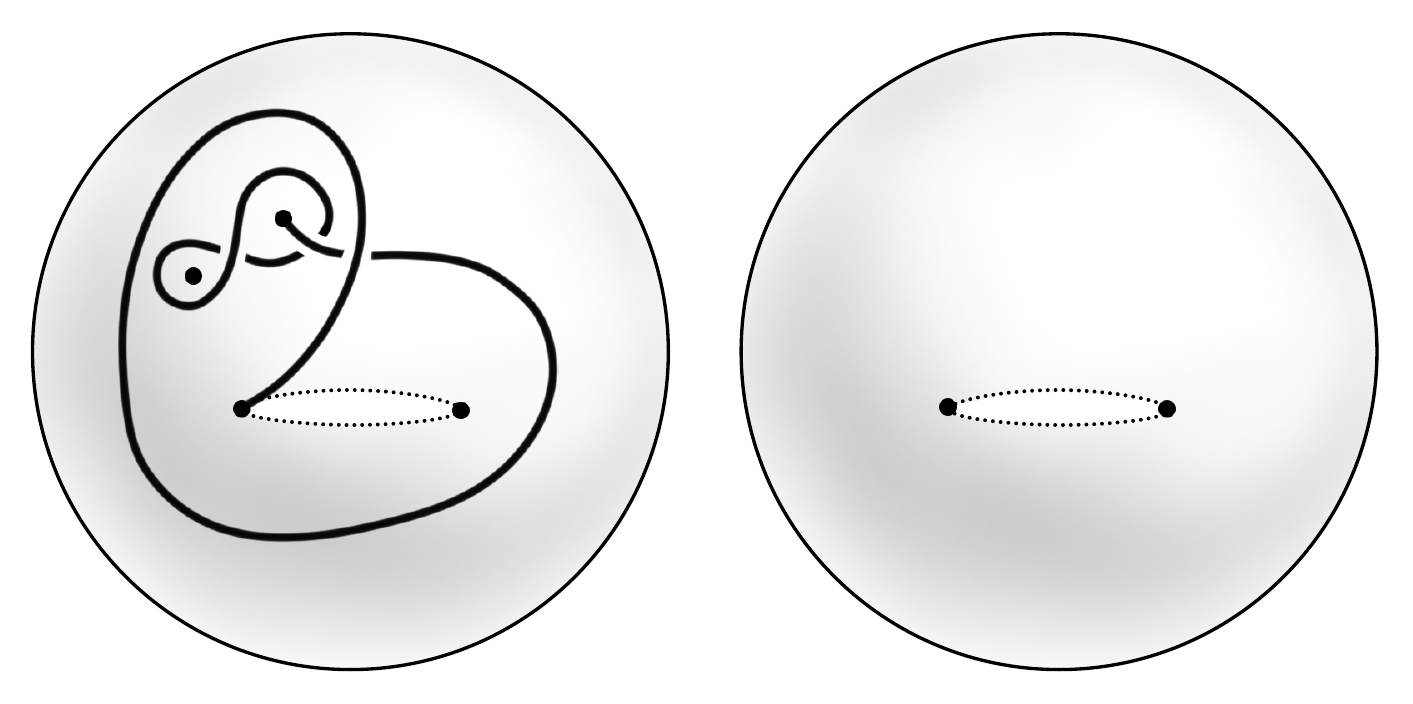}
        \caption{The diagram $\diag_0'$ for the diagram $\diag_0$ from Figure \ref{diagZero}, shown in the same schematic as Figure \ref{doubleCover}. The diagrams $\diag_0$ and $\diag_0'$ are ambiently isotopic.}
        \label{doublecoverdiagzero}
    \end{figure}
\end{proof}

We conclude this section with some open questions.

\begin{question}\label{all-heights}
    Does every generalized knotoid $\kappa$ have a diagram $\diag$ that simultaneously realizes the height between every pair of poles of $\kappa$?
\end{question}

As a corollary to Theorem \ref{heightTheorem}, we obtain an affirmative answer to Question \ref{all-heights} in the case that the height spectrum of $\kappa$ is the zero matrix: we may take $\diag$ to be a minimal-crossing diagram of $\kappa$.

\begin{question}\label{minCrossingArbitraryHeight}
    Let $\kappa$ be a generalized knotoid and let $A,B \in P(\kappa)$ be poles. Does every minimal-crossing diagram of $\kappa$ realize the height between $A$ and $B$?
\end{question}

Question \ref{minCrossingArbitraryHeight} asks a generalization of Theorem \ref{heightTheorem} from height $0$ to arbitrary heights, and an affirmative answer to Question \ref{minCrossingArbitraryHeight} implies an affirmative answer to Question \ref{all-heights}.

\begin{question}
    Does Theorem \ref{heightTheorem} generalize to non-classical generalized knotoids (for example, generalized knotoids on a genus $g$ surface)?
\end{question}

The proof we have given for Theorem \ref{heightTheorem} does not immediately generalize to higher genus surfaces because it uses the fact that there is a double branched cover of the sphere by itself. Generalizing this method may require rephrasing the argument in terms of a ``virtual theory'' for generalized knotoids.

\begin{question}
    Which matrices of integers are realized as the height spectrum of a generalized knotoid?
\end{question}

\section{Index polynomials for Generalized Knotoids}\label{index}
\subsection{Notation}
\begin{figure}
    \captionsetup[subfigure]{labelformat=empty}
    \centering
    \begin{subfigure}{0.3\textwidth}
    \centering
    \includegraphics[width=2cm]{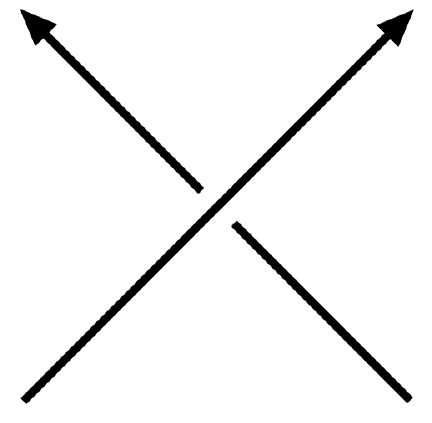}
    \caption{$\on{sgn}(c) = 1$}
    \end{subfigure}
    \begin{subfigure}{0.3\textwidth}
    \centering
    \includegraphics[width=2cm]{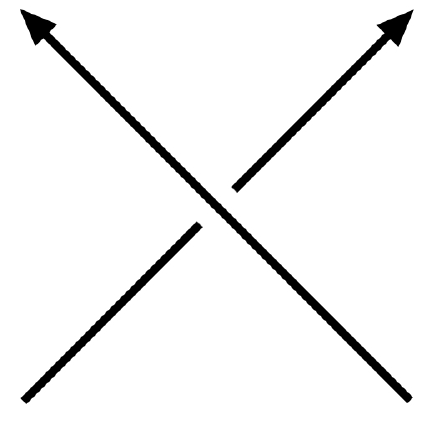}
    \caption{$\on{sgn}(c) = -1$}
    \end{subfigure}
    \caption{Positive and negative crossing, respectively}
    \label{crossing-sign}
\end{figure}

\begin{figure}
    \centering
    \begin{subfigure}{0.3\textwidth}
    \centering
    \includegraphics[width=2cm]{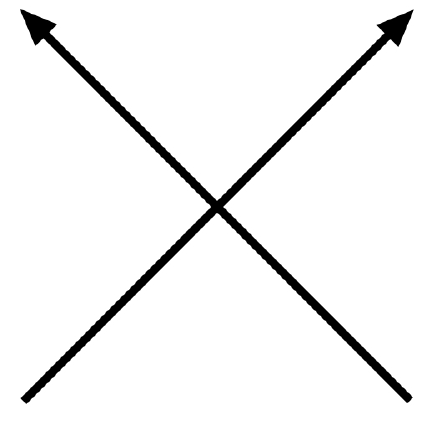}
    \end{subfigure}
    \begin{subfigure}{0.3\textwidth}
    \centering
    \includegraphics[width=2cm]{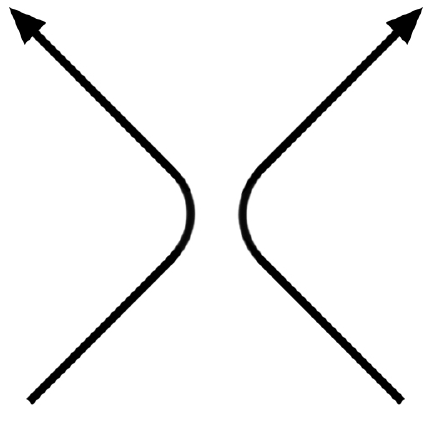}
    \end{subfigure}
    \caption{A crossing and its oriented smoothing. The smoothing does not depend on over/under data.}
    \label{oriented-smoothing}
\end{figure}

Given a crossing $c$ in an oriented diagram, we let $\on{sgn}(c) \in \{1, -1\}$ denote the \emph{sign} of the crossing. See Figure \ref{crossing-sign}.
The \emph{oriented smoothing} at $c$ is the smoothing that respects the orientation of the incoming and outgoing strands. See Figure \ref{oriented-smoothing}.
Given a diagram $\diag$ and two constituents $\alpha, \beta \in C(\diag)$, we define their \emph{linking number} by the half-integer
\[\on{lk}(\alpha, \beta) = \frac 12\sum_c \on{sgn}(c),\]
where the sum is taken over all crossings between $\alpha$ and $\beta$. Observe that the linking number between $\alpha$ and $\beta$ is an invariant of the generalized knotoid represented by $\diag$.

Given two (possibly identical) oriented generalized knotoid diagrams $\diag_1, \diag_2$ considered on the same surface $\Sigma$ and constituents $\alpha \in C(\diag_1)$ and $\beta \in C(\diag_2)$, let $T(\alpha, \beta)$ denote the finite subset of $\alpha \cap \beta$ consisting of transverse intersections between $\alpha$ and $\beta$, excluding potential intersections at poles. At each $c \in T(\alpha, \beta)$, the tangent vectors $v_\alpha, v_\beta$ of $\alpha$ and $\beta$ form a basis for the tangent space of $\Sigma$ at $c$. We set $\varepsilon(c) = 1$ if $(v_\beta, v_\alpha)$ is a positively oriented basis, and $\varepsilon(c) = -1$ otherwise. (Equivalently, $\varepsilon(c)$ gives the sign of $c$ if $\beta$ is taken to be the over-strand.) The \emph{algebraic intersection number} of $\alpha$ and $\beta$, denoted $\alpha \cdot \beta$, is defined as
\[\alpha \cdot \beta = \sum_{c \in T(\alpha, \beta)} \varepsilon(c).\]
The algebraic intersection number does not require over/under data. The sign convention we have adopted matches the convention of the algebraic intersection number defined in \cite{Turaev1}.

Proposition \ref{loops-on-sphere} and Corollary \ref{path-independence} are well known.
\begin{proposition}\label{loops-on-sphere}
If $\alpha$ and $\beta$ are loop constituents on $S^2$ such that all intersections in $\alpha \cap \beta$ are transverse and away from crossings, then $\alpha \cdot \beta = 0$.
\end{proposition}
\begin{corollary}\label{path-independence}
Let $\beta$ be a loop constituent on $S^2$ and fix two points $A,B$ in $S^2 \setminus \beta$. Let $\alpha$ be a segment constituent between $A$ and $B$ intersecting $\beta$ transversely and away from crossings. Then $\alpha \cdot \beta$ does not depend on the choice of $\alpha$.
\end{corollary}
%\begin{proof}
%    Let $\beta_1, \beta_2$ be two such choices of $\beta$. Let $\gamma$ denote the loop formed by concatenating $\beta_1$ with orientation-reversed $\beta_2$. By Proposition \ref{loops-on-sphere}, we have
%    \begin{align*}
%        0 &= \gamma \cdot \alpha\\
%        &= \beta_1 \cdot \alpha - \beta_2 \cdot \alpha.
%    \end{align*}
%\end{proof}

The \emph{index polynomial} is an invariant for virtual knots introduced by Turaev in \cite{Turaev2} and Henrich in \cite{henrich}. A similar \emph{affine index polynomial} for oriented virtual knots was defined by Kauffman in \cite{affine-index-virtual} and extended to knotoids in \cite{NIV}. Kim, Im, and Lee defined an index polynomial for knotoids in \cite{family-index-poly-knotoid} that is distinct from the affine index polynomial.

These polynomial invariants are based on the following idea: consider a classical crossing $c$ in a diagram of a knotoid or virtual knot and apply an oriented smoothing at $c$. The smoothing produces a diagram with two constituents $\alpha_c, \beta_c$. In the case of a virtual knot, the smoothed diagram is a two-component virtual link, and in the case of a knotoid, it is a loop constituent and a segment constituent. An \emph{intersection index}, denoted $\on{in}(c)$, is computed from the algebraic intersection number of $\alpha_c$ and $\beta_c$, and the polynomial invariant takes the form
\begin{equation}\label{base-index-poly}F_\diag(t) = \sum_{c} \on{sgn}(c)(t^{\on{in}(c)} - 1),\end{equation}
where the sum is taken over all crossings of the diagram $\diag$. The values of $\on{in}(c)$ computed for the various index polynomials differ only by sign:
\begin{enumerate}[(i)]
    \item In Henrich's index polynomial for virtual knots, the intersection index is taken to be $\on{in}(c) = |\alpha_c \cdot \beta_c|$. 
    \item In Kauffman's affine index polynomial $F^{\on{aff}}(t)$ for oriented virtual knots and knotoids, the constituents $\alpha_c, \beta_c$ are labeled such that $\alpha_c$ contains the incoming portion of the over-strand and $\beta_c$ contains the incoming portion of the under-strand. The intersection index is taken to be $\on{in}(c) = w(c) := \alpha_c \cdot \beta_c$. (This is not Kauffman's original formulation of the affine index polynomial in \cite{affine-index-virtual}, but it was shown to be equivalent in \cite{wriggle}.)
    \item In Kim, Im, and Lee's index polynomial $F^{\on{ind}}(t)$ for knotoids, the constituents $\alpha_c, \beta_c$ are labeled so that $\alpha_c$ is the loop constituent and $\beta_c$ is the segment constituent, and the intersection index is taken to be $\on{in}(c) = \on{ind}(c) := \alpha_c \cdot \beta_c$. See the discussion in Section 5 of \cite{linov}.
\end{enumerate}
In \cite{NIV}, it is shown that for a knotoid $k$, the affine index polynomial $F^{\on{aff}}_k(t)$ satisfies $F^{\on{aff}}_k(t) = F^{\on{aff}}_k(t^{-1})$. It follows that the affine index polynomial for knotoids is recoverable from the index polynomial via
\[F^{\on{aff}}_k(t) = \frac 12 (F^{\on{ind}}_k(t) + F^{\on{ind}}_k(t^{-1})).\]
In \cite{family-index-poly-knotoid}, the index polynomial is strengthened to a two-variable polynomial
\[F^{\on{ind}}_\diag(s,t) = \sum_{c \in U(\diag)} \on{sgn}(c)(t^{\on{ind}(c)} -1) + \sum_{c \in O(\diag)} \on{sgn}(c)(s^{\on{ind}(c)} - 1).\]
Here, $U(\diag)$ (resp. $O(\diag)$) denotes the set of \emph{early undercrossings} (resp. \emph{early overcrossings}) of the diagram $\diag$, the crossings $c$ such that $c$ is first encountered as an undercrossing (resp. overcrossing).

\subsection{Generalized index polynomials}\label{genindexpoly}
We now define index polynomial invariants for generalized knotoids. Let $\diag$ be an oriented, constituent-labeled generalized knotoid diagram. Let $c$ be a crossing, and let $o_c$ and $u_c$ denote the constituents of $\diag$ containing the over-strand and under-strand at $c$, respectively. Suppose an oriented smoothing is applied to $c$. Let $\alpha_c$ denote the constituent of the smoothed diagram containing the incoming portion of $o_c$, and let $\beta_c$ denote the constituent of the smoothed diagram containing the incoming portion of $u_c$. Note that $o_c$ and $u_c$ are not necessarily distinct, and $\alpha_c$ and $\beta_c$ are not necessarily distinct.

\begin{figure}[htbp]
    \captionsetup[subfigure]{labelfont=rm}
    \begin{subfigure}{0.48\textwidth}
    \centering
    \includegraphics[width=5.75cm]{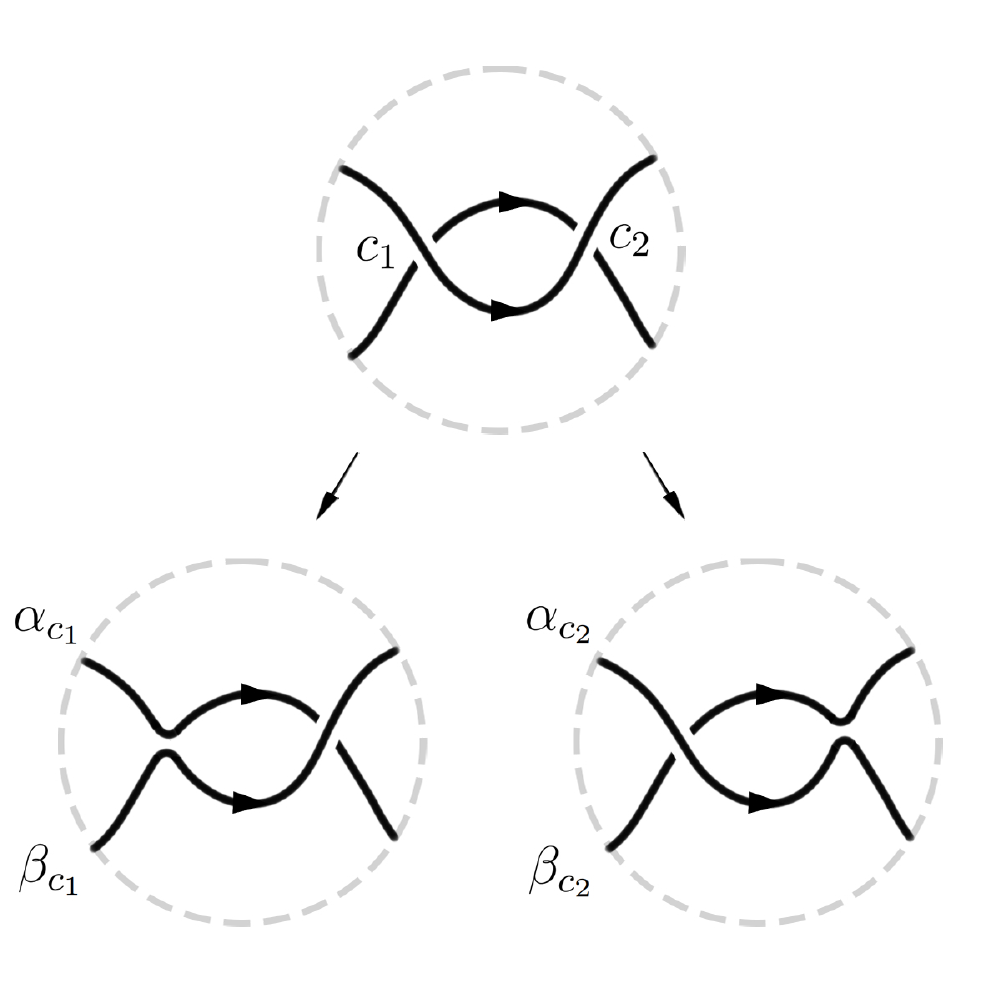}
    \caption{}
    \label{r2invariancea}
    \end{subfigure}
    \begin{subfigure}{0.48\textwidth}
    \centering
    \includegraphics[width=5.75cm]{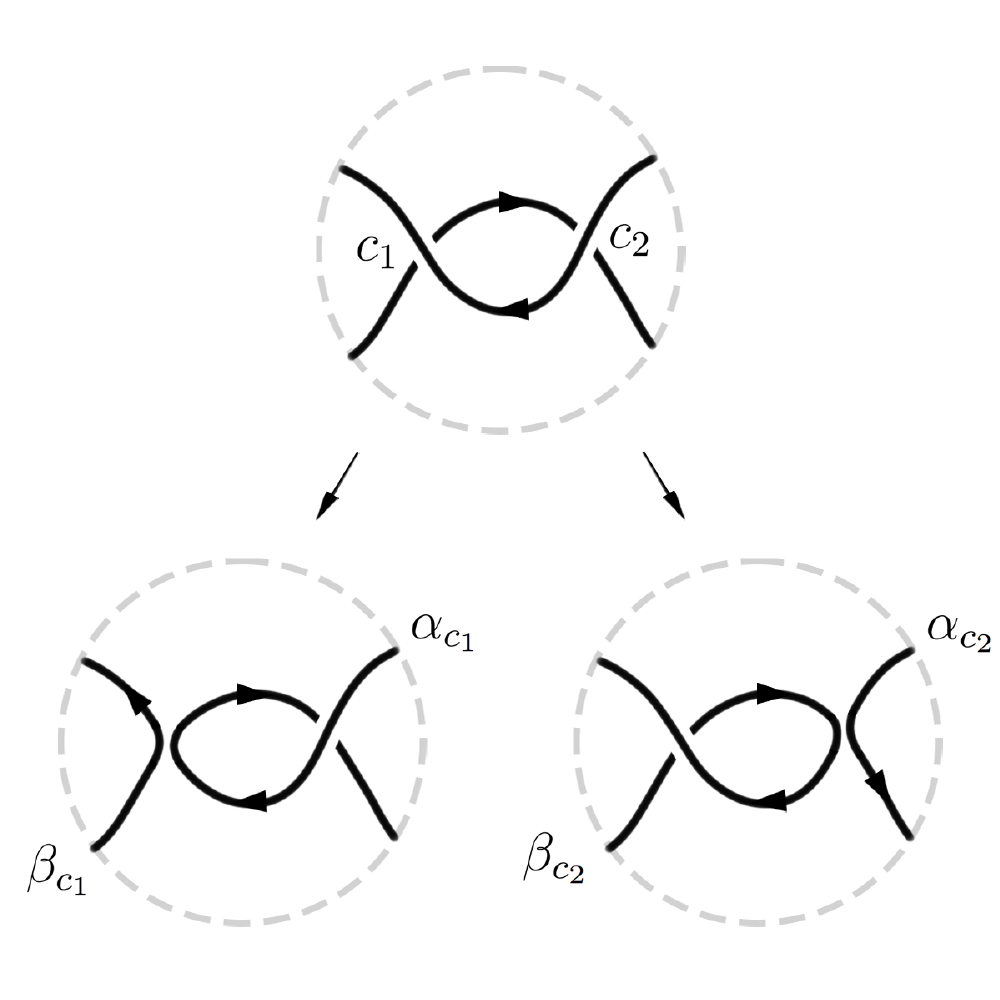}
    \caption{}
    \label{r2invarianceb}
    \end{subfigure}
    
    \caption{Two possible configurations for a Type II Reidemeister move. In each of Figures \ref{r2invariancea} and \ref{r2invarianceb}, the two diagrams obtained by smoothing at $c_1$ and $c_2$ carry the same algebraic intersection data.}
    \label{r2invariance}
\end{figure}

We introduce variables $r_{e_1, e_2}$, $t_{e_1}$, and $s_{e_1}$ for every pair of (not necessarily distinct) constituents $e_1, e_2 \in C(\diag)$ and define
\[g(c) = r_{o_c, u_c}\left(\prod_{e \in C(\diag)} t_e^{\alpha_c \cdot e} - 1\right)\left(\prod_{e \in C(\diag)}s_e^{\beta_c \cdot e} - 1\right).\]
We define the \emph{generalized index polynomial} by
\begin{align}\label{defgip}
G_\diag(r, s, t) = \sum_c \on{sgn}(c)g(c).\end{align}
Here $r$ represents a vector whose components are the variables $r_{e_1, e_2}$ in some fixed order, and $s$ and $t$ are similar. Operations on $r$, $s$, and $t$ are applied component-wise. For example, the polynomial $G_\diag(r, s^{-1}, t)$ is obtained by replacing $s_e$ with $s_e^{-1}$ for all $e \in C(\diag)$, and $G_\diag(r, t, s)$ is obtained by swapping $t_e$ and $s_e$ for all $e \in C(\diag)$. For $r$, we let $\overline r$ indicate the operation that swaps $r_{e_1, e_2}$ and $r_{e_2, e_1}$ for all constituents $e_1, e_2$.
\begin{theorem}\label{genindexinvariance}
    The polynomial $G_\diag$ is a generalized knotoid invariant.
\end{theorem}
\begin{proof}
    Suppose a generalized knotoid $\kappa$ has diagrams $\diag, \diag'$ related by a Reidemeister move. The crossings of $\diag$ and $\diag'$ are in natural correspondence except potentially for crossings created or destroyed by a Type I or II Reidemeister move. For corresponding crossings, the values of $\on{sgn}(c)$, $o_c$, $u_c$, $\alpha_c \cdot e$, and $\beta_c \cdot e$ are identical for any $e \in C(\kappa)$. It remains to verify that the crossings created or destroyed by a Type I or Type II Reidemeister move do not contribute to $G_\diag$. 
    \begin{enumerate}[(i)]
        \item Consider a Type I Reidemeister move, and assume without loss of generality that a crossing $c$ is created. Observe that at least one of $\alpha_c, \beta_c$ is a loop with no crossings, so $g(c) = 0$.
        \item Consider a Type II Reidemeister move, and assume without loss of generality that two crossings $c_1, c_2$ are created. Observe that $c_1$ and $c_2$ have opposite sign and satisfy $(o_{c_1}, u_{c_1}) = (o_{c_2}, u_{c_2})$. Outside of the neighborhood where the Reidemeister move is performed, the constituents $\alpha_{c_1}$ and $\alpha_{c_2}$ are identical, and the constituents $\beta_{c_1}$ and $\beta_{c_2}$ are identical. Inside the neighborhood, the oriented smoothings at $c_1$ and $c_2$ each leave one crossing, the two of which are indistinguishable by algebraic intersection numbers. See Figure \ref{r2invariance}. Thus $\alpha_{c_1} \cdot e = \alpha_{c_2} \cdot e$ and $\beta_{c_1} \cdot e = \beta_{c_2} \cdot e$ for any $e$, and it follows that $g(c_1)= g(c_2)$. Finally, $\on{sgn}(c_1)g(c_1) + \on{sgn}(c_2)g(c_2) = 0$, as desired.
    \end{enumerate}
\end{proof}

The generalized index polynomial extends the affine index polynomial in the following sense. 
\begin{proposition}\label{reduce-to-aip}
    Two knotoids have identical affine index polynomials if and only if they have identical generalized index polynomials.
\end{proposition}
\begin{proof}
    If $\diag$ is a spherical or planar knotoid diagram, the generalized index polynomial becomes a three-variable polynomial
    \[G_\diag(r,s,t) = \sum_c \on{sgn}(c)r(t^{\alpha_c \cdot e} - 1)(s^{\beta_c \cdot e} -1 ).\]
    Here, $e$ is the single segment constituent of $\diag$. Observe that $\alpha_c \cdot e = \alpha_c \cdot \beta_c = -\beta_c \cdot e = w(c)$. It follows that the affine index polynomial $F^{\on{aff}}_\diag(t)$ is recovered from $G_\diag(s,t)$ by
    \[F^{\on{aff}}_\diag(t) = -G_\diag(1, 0, t).\]
    Conversely, $G_\diag(r,s,t)$ is recovered from $F^{\on{aff}}_\diag(t)$ by
    \[G_\diag(r,s,t) = r(F^{\on{aff}}_\diag(ts^{-1}) - F^{\on{aff}}_\diag(t) - F^{\on{aff}}_\diag(s^{-1})).\]
\end{proof}
By Proposition \ref{loops-on-sphere}, the generalized index polynomial is trivial for any classical link.

If a generalized knotoid diagram $\diag$ is the disjoint union of diagrams $\diag_1$ and $\diag_2$ with no crossings occuring between constituents of $\diag_1$ and $\diag_2$, then the generalized index polynomial $G_\diag$ provides no information on their relative positioning. For example, a valency-zero pole can be moved to any region of a diagram without changing the generalized index polynomial.

We extend the generalized index polynomial for pole-labeled generalized knotoids with at least one pole. Given a diagram $\diag$ and a fixed pole $P \in P(\diag)$, we choose oriented shortcuts $\{\gamma_Q \mid Q \in P(\diag)\}$ from $P$ to each pole in $\diag$. For a constituent $e$, let $\mathbf 1_L(e)$ denote the indicator function that is $1$ if $e$ is a loop constituent and $0$ otherwise. We introduce variables $a_Q, b_Q$ for each $Q \in P(\diag)$ and define, for each crossing $c$, the polynomial
\[h_P(c) = \left(\prod_{Q \in P(\diag)} a_Q^{\gamma_Q \cdot \alpha_c}\right)^{\mathbf 1_L(\alpha_c)}\left(\prod_{Q \in P(\diag)} b_Q^{\gamma_Q \cdot \beta_c}\right)^{\mathbf 1_L(\beta_c)}.
\]
We define the \emph{base-pointed index polynomial} by
\begin{equation}\label{bpgip}
    G_{\diag, P}(r,s,t,a,b) = \sum_c \on{sgn}(c)g(c)h_P(c).
\end{equation}
The notational remarks following Equation (\ref{defgip}) apply to the variables $a$ and $b$ as well.

By Corollary \ref{path-independence}, the base-pointed index polynomial is independent of the choice of shortcuts $\gamma_Q$, and the proof that $G_{\diag, P}$ is a generalized knotoid invariant is identical to the proof of Theorem \ref{genindexinvariance}. The choice of the base-point pole $P$ is immaterial. If $P'$ is another pole, the polynomial $G_{\diag, P'}$ is obtained from $G_{\diag, P}$ by the following procedure: modify each monomial in $G_{\diag, P}$ by subtracting the exponent of $a_{P'}$ from the exponent of each $a_Q$ and subtracting the exponent of $b_{P'}$ from the exponent of each $b_Q$.

\begin{figure}[htbp]
    \captionsetup[subfigure]{labelfont=rm}
    \begin{subfigure}{0.3\textwidth}
    \centering
    \includegraphics[width=3.3cm]{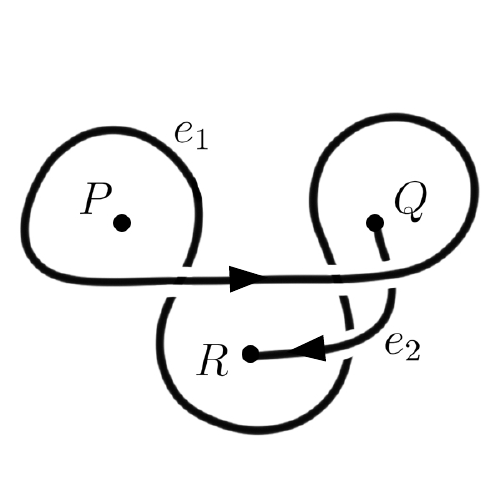}
    \caption{}
    \label{mickey}
    \end{subfigure}
    \begin{subfigure}{0.3\textwidth}
    \centering
    \includegraphics[width=3.3cm]{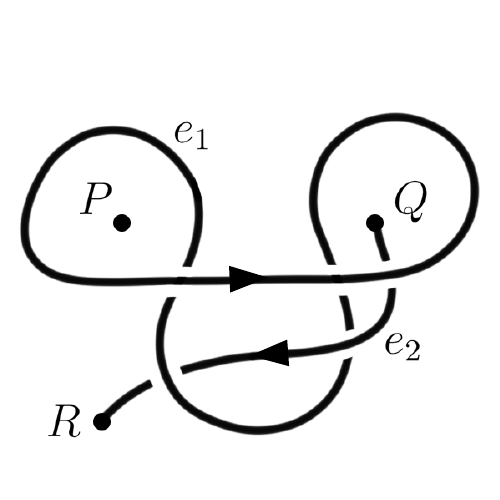}
    \caption{}
    \label{mickeytwo}
    \end{subfigure}
    \begin{subfigure}{0.3\textwidth}
    \centering
    \includegraphics[width=3.3cm]{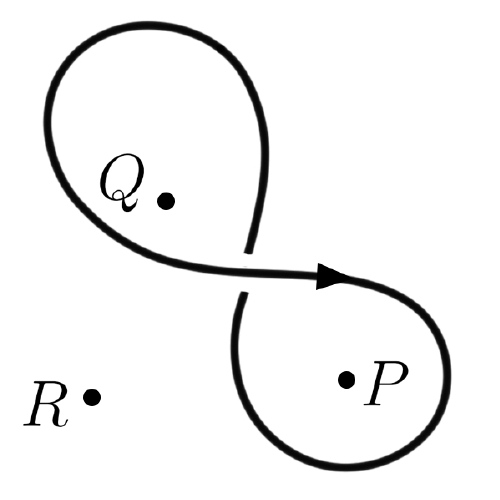}
    \caption{}
    \label{infty1}
    \end{subfigure}
    \begin{subfigure}{0.3\textwidth}
    \centering
    \includegraphics[width=3.3cm]{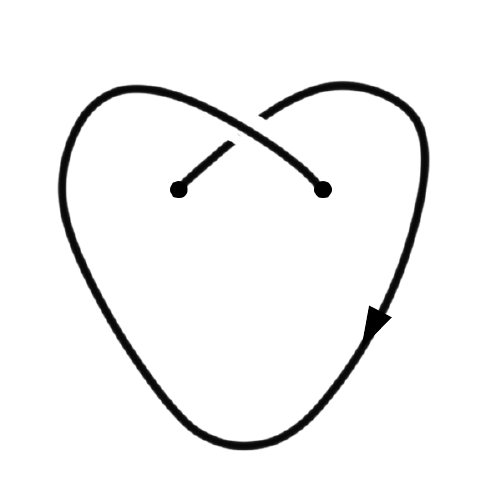}
    \caption{}
    \label{nontrivialplanar}
    \end{subfigure}
    \begin{subfigure}{0.3\textwidth}
    \centering
    \includegraphics[width=3.3cm]{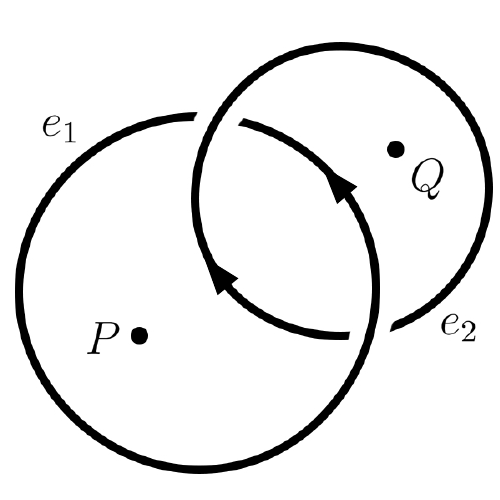}
    \caption{}
    \label{stakedhopf}
    \end{subfigure}
    \caption{}
    \label{indexpolyexamples}
\end{figure}

\begin{example}\label{mickeyexample}
    The generalized knotoid $\kappa$ in Figure \ref{mickey} has
    \[G_{\kappa, P} = r_{e_1, e_1}((s_{e_2}^2 - 1)a_Q^{-1}a_R^{-1}b_Qb_R^{-1} - (t_{e_2} - 1)(s_{e_2}-1)a_Q^{-1}a_R^{-2}b_Q)\]
    and
    \[G_{\kappa, R} = r_{e_1, e_1}((s_{e_2}^2 - 1)a_Pb_Q^2 - (t_{e_2} - 1)(s_{e_2}-1)a_P^2a_Q^1b_Q).\]
\end{example}

\begin{example}
    The generalized knotoid $\kappa$ in Figure \ref{mickeytwo} has $G_{\kappa, P} = 0$, so the base-pointed generalized index polynomial cannot distinguish $\kappa$ from a diagram with no crossings. However, every diagram of $\kappa$ has at least one crossing since $\on{lk}(e_1, e_2) = 1/2$.
\end{example}
We now define a variant of the base-pointed index polynomial. With all notation as before, define
\[\tilde g(c) = r_{o_c, u_c}\prod\limits_{e \in C(\diag)}t_e^{\alpha_c \cdot e}s_e^{\beta_c \cdot e}\]
and
\[\tilde h_P(c) = \left(\mathbf 1_L(\alpha_c)\prod_{Q \in P(\diag)} a_Q^{\gamma_Q \cdot \alpha_c} - 1\right)\left(\mathbf 1_L(\beta_c)\prod_{Q \in P(\diag)} b_Q^{\gamma_Q \cdot \beta_c} -1\right).\]
We define the \emph{pole-centric base-pointed index polynomial} by
\[\tilde G_{\diag, P}(r,s,t,a,b) = \sum_c \on{sgn}(c)\tilde g(c) \tilde h_P(c).\]
The proof that $\tilde G_{\diag, P}$ is a generalized knotoid invariant is nearly identical to the proof of Theorem \ref{genindexinvariance}. The remarks following Equation~\eqref{bpgip} for $G_{\diag, P}$ hold for $\tilde G_{\diag, P}$ as well.

In analogy to Proposition \ref{reduce-to-aip}, the pole-centric base-pointed index polynomial extends the strengthened index polynomial $F^{\on{ind}}(s,t)$ for knotoids. Suppose $k_1$ and $k_2$ are two spherical knotoids, each with poles labeled $L$ and $H$ and a constituent oriented from $L$ to $H$.

\begin{proposition}
     If $G_{k_1, L} = G_{k_2, L}$, then $F^{\on{ind}}_{k_1}(s,t) = F^{\on{ind}}_{k_1}(s,t)$. 
\end{proposition}
\begin{proof}
For a knotoid diagram $\diag$ oriented from $L$ to $H$, the pole-centric base-pointed index polynomial based at $L$ is a five-variable polynomial $\tilde G_{\diag, L}(r, s, t, a_H, b_H)$. Set $r=s=t = 1$ to obtain a two-variable polynomial $\tilde G_{\diag, L}(a_H, b_H)$. For each crossing $c$ in $\diag$, exactly one of $\alpha_c, \beta_c$ is a loop constituent. The crossings for which $\alpha_c$ is the loop constituent is precisely the set of early undercrossings, and the crossings for which $\beta_c$ is the loop constituent is precisely the set of early overcrossings. Thus
\[\tilde G_{\diag, L}(a_H, b_H) = -\sum_{c \in U(k)} \on{sgn}(c)(a_H^{\gamma_{H} \cdot \alpha_c} - 1)  - \sum_{c \in O(\diag)} \on{sgn}(c)(b_H^{\gamma_{H} \cdot e} - 1).\]
Suppose $\alpha_c$ is the loop constituent. Observe that $\beta_c$ defines an oriented shortcut from $L$ to $H$. It follows from Corollary \ref{path-independence} that
\[\gamma_H \cdot \alpha_c = \beta_c \cdot \alpha_c = -\on{ind(c)}.\]
Similarly, if $\beta_c$ is the loop constituent, then $\gamma_H \cdot \beta_c = \alpha_c \cdot \beta_c = -\on{ind}(c)$. Thus we obtain
\begin{align*}
    -\tilde G_{\diag, L}(s^{-1}, t^{-1}) &= \sum_{c \in U(\diag)} \on{sgn}(c)(s^{\on{ind}(c)} - 1) + \sum_{c \in O(\diag)}\on{sgn}(c)(t^{\on{ind}(c)} - 1)\\
    &= F^{\on{ind}}_\diag(s,t).
\end{align*}
In particular, the strengthened index polynomial is recoverable from the pole-centric base-pointed index polynomial, and the conclusion follows.
\end{proof}

\begin{example}
    The staked knot $\kappa$ from Figure \ref{infty1} has $\tilde G_{\kappa, P} = a_Qb_Qb_R$. On the other hand, the base-pointed index polynomial $G_{\kappa, P}$ is trivial for any staked link $\kappa$ by Proposition \ref{loops-on-sphere}.
\end{example}
\begin{question}
    Does there exist a pole-labeled, constituent-labeled, oriented generalized knotoid $\kappa$ with a pole $P$ such that $\tilde G_{\kappa, P} = 0$ but $G_{\kappa, P} \neq 0$?
\end{question}

The index polynomials of the form given by Equation~\eqref{base-index-poly} cannot distinguish between inequivalent planar knotoids with diagrams that are equivalent when considered on the sphere. By viewing a planar knotoid as a generalized knotoid with three poles $L$, $H$, and $\infty$ and a constituent oriented from $L$ to $H$, as in Example \ref{planar_example}, the base-pointed index polynomial and its pole-centric variant yield new invariants for planar knotoids. For example, the planar knotoid $k$ in Figure \ref{nontrivialplanar} has trivial index polynomial but has $\tilde G_{k, \infty} = a_La_H$.

Following the notation of \cite{Turaev1} and \cite{linov}, we extend the basic knotoid involutions to generalized knotoids and state their effects on the generalized index polynomials. Given a diagram $\diag$, we let $\on{mir}(\diag)$ denote the diagram obtained by toggling all over/under crossing data, and we let $\on{sym}(\diag)$ denote the diagram obtained by reflecting across a great circle of $S^2$ (and preserving crossing data). We let $\on{rot}$ denote the composition $\on{mir} \circ \on{sym} = \on{sym} \circ \on{mir}$. For an oriented diagram $\diag$, we let $\on{rev}(\diag)$ denote the diagram obtained by reversing the orientation on all constituents. It is clear that $\on{mir}$, $\on{sym}$, $\on{rot}$, and $\on{rev}$ define involutions on the set of generalized knotoids, possibly labeled or oriented. Note that $\on{mir}$, $\on{sym}$, and $\on{rot}$ preserve the underlying directed graph of an oriented generalized knotoid, while $\on{rev}$ may not.
\begin{proposition}
Let $\kappa$ be pole-labeled, constituent-labeled, oriented generalized knotoid and let $P \in P(\kappa)$. The base-pointed index polynomial satisfies:
\begin{enumerate}[(a)]
    \item $G_{\on{mir}(\kappa), \on{mir}(P)}(r,s,t,a,b) = -G_{\kappa, P}(\overline r, t, s , b, a)$,
    \item $G_{\on{sym}(\kappa), \on{sym}(P)}(r,s,t,a,b) = -G_{\kappa, P
}(r, s^{-1}, t^{-1}, a^{-1}, b^{-1})$,
    \item $G_{\on{rot}(\kappa), \on{rot}(P)}(r,s,t,a,b) = G_{\kappa, P}(\overline r, t^{-1}, s^{-1}, b^{-1}, a^{-1})$.
\end{enumerate}
Identical relationships hold for the pole-centric variant.
\end{proposition}

\begin{proof}
    Part (a) follows from the fact that, for each crossing $c$, the operation $\on{mir}$ negates $\on{sgn}(c)$ and swaps $\alpha_c$ and $\beta_c$. Part (b) follows from the fact that $\on{sym}$ negates both the sign of each crossing and all algebraic intersection numbers. Part (c) follows from (a) and (b).
\end{proof}
Thus the base-pointed index polynomials may be used to distinguish between $\kappa$ and its images under the basic involutions. 

The base-pointed index polynomials also yield lower bounds on the height between poles in a generalized knotoid. Given a Laurent polynomial $F \in \mathbb Z[x_1^{\pm 1}, \dots, x_n^{\pm 1}]$, we let $[x_i]F \in \mathbb Z[x_1^{\pm 1}, \dots, x_{i-1}^{\pm 1}, x_{i+1}^{\pm 1}, \dots, x_n^{\pm 1}]$ denote the coefficient of $x_i$ in $F$. For any subset $S = \{x_{i_1}, \dots, x_{i_m}\}$ of the variables and a monomial $M = \prod_{j=1}^n x_j^{u_j}$, we define $\deg_S(M)=\sum_{k=1}^m |u_{i_k}|$, and we let $\deg_S(F)$ denote the maximum of $\deg_S(M)$ over all monomials $M$ with nonzero coefficient $F$. (We omit the set notation from $S$ for brevity).
\begin{proposition}\label{heightbound} %and I'm homebound~~
     Let $\kappa$ be a pole-labeled, constituent-labeled, oriented generalized knotoid and let $P, Q \in P(\kappa)$ be poles. Then
     \[h_\kappa(P,Q) \ge \sum_{e \in C(\kappa)} \deg_{a_Q, b_Q}([r_{e, e}]G_{\kappa, P}).\]
     An identical inequality holds with $G_{\kappa, P}$ replaced by $\tilde G_{\kappa, P}$.
\end{proposition}
\begin{proof}
    The following argument holds for both $G_{\kappa, P}$ and $\tilde G_{\kappa, P}$. Let $\diag$ be any diagram representing $\kappa$ and let $\gamma$ be any shortcut oriented from $P$ to $Q$. It suffices to show that, for any constituent $e \in C(\diag)$, the number of intersections $\#(\gamma \cap e)$ of $\gamma$ with $e$ away from poles is at least $d := \deg_{a_Q, b_Q}([r_{e,e}]G_{\kappa, P})$.
    
    If $d = 0$, this is trivial. Assuming $d > 0$, the diagram $\diag$ must contain a crossing $c$ with $o_c = u_c = e$ and $\mathbf 1_L(\alpha_c)|\gamma_Q \cdot \alpha_c| + \mathbf 1_L(\beta_c)|\gamma_Q \cdot \beta_c| = d$. Note that $\alpha_c \neq \beta_c$, and
    \begin{align*}
        \#(\gamma \cap e) &= \#(\gamma \cap \alpha_c) + \#(\gamma \cap \beta_c)\\
        &\ge |\gamma_Q \cdot \alpha_c| + |\gamma_Q \cdot \beta_c|\\
        &\ge d.
    \end{align*}
\end{proof}
Proposition \ref{heightbound} implies that the generalized knotoid $\kappa_1$ from Example \ref{mickeyexample} has $h_{\kappa_1}(P, Q) = h_{\kappa_1}(P, R) = h_{\kappa_1}(R,Q) = 2$, and the generalized knotoid $\kappa_2$ from Example \ref{infty1} has $h_{\kappa_2}(P, Q) = 2$ and $h_{\kappa_2}(P, R) = h_{\kappa_2}(Q,R) = 1$.

For certain generalized knotoids $\kappa$, the polynomials $G_{\kappa, P}$ or $\tilde G_{\kappa, P}$ give more information than the bound stated in Proposition \ref{heightbound}, as the bound only considers intersections of a shortcut $\gamma$ with loops $\alpha_c, \beta_c$ that result from smoothing at a self-intersection of a constituent $e$. The generalized knotoid $\kappa$ in Figure \ref{stakedhopf} has $\tilde G_{\kappa, P} = 2r_{e_1, e_2}a_Q^{-2}$, and an argument similar to the proof of Proposition \ref{heightbound} shows that $h_\kappa(P,Q) = 2$ despite the fact that Proposition \ref{heightbound} only yields the trivial bound $h_\kappa(P, Q) \ge 0$.

We remark on non-constituent-labeled generalized knotoids. If $\kappa$ is simply an oriented, generalized knotoid, we obtain a two-variable generalized index polynomial $G_{\kappa}(s,t)$ by computing the generalized index polynomial of $\kappa$ with constituents arbitrarily labeled, and then setting $r_{e_1, e_2}=1$, $t_{e_1}=t$, and $s_{e_1} = s$ for all $e_1, e_2 \in C(\kappa)$. In the case that $\kappa$ is pole-labeled, we again label the constituents arbitrarily and compute the base-pointed index polynomial (or the pole-centric variant). We then re-index the variables $r_{e_1, e_2}, t_e, s_e$ as follows. Each subscript $e$ corresponding to a segment constituent $e$ oriented from pole $P$ to pole $Q$ is replaced by $(P, Q)$, and all subscripts corresponding to loop constituents are replaced by a single symbol $L$.

\section{A Bracket Polynomial for Generalized Knotoids}\label{bracket}
\subsection{Generalized bracket polynomial}
We define an extension of the Kauffman bracket polynomial for oriented, pole-labeled, edge-labeled generalized knotoids.

\begin{figure}
    \centering
    \includegraphics[width=2cm]{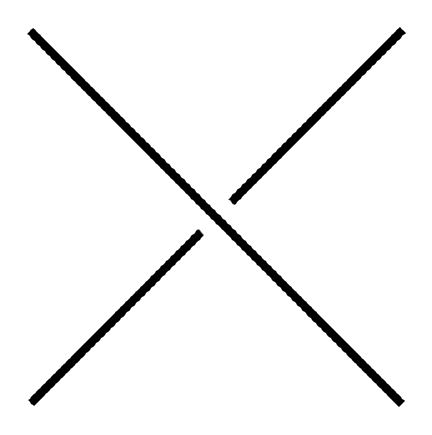}
    \includegraphics[width=2cm]{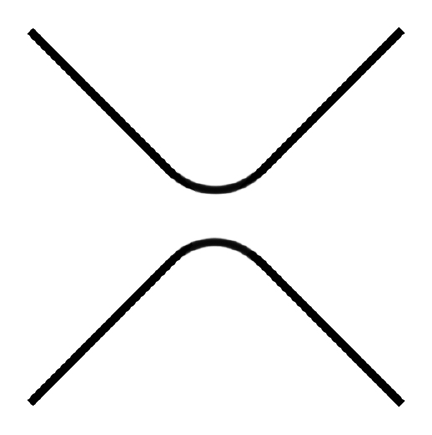}
    \includegraphics[width=2cm]{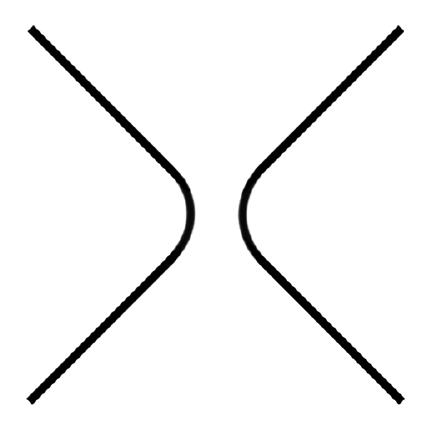}
    \caption{A crossing, its $A$-type smoothing, and its $B$-type smoothing, respectively.}
    \label{ABsmoothings}
\end{figure}

Given a generalized knotoid diagram $\diag$ and a crossing $c$ in $\diag$, a \emph{smoothing at $c$} is one of the two operations pictured in Figure \ref{ABsmoothings}, and is either $A$-type or $B$-type.
Observe that a smoothing preserves the number of segment constituents. A \emph{state} $s$ of $\diag$ is a diagram obtained by smoothing at each of the crossings of $\diag$. Let $\sigma(s)$ denote the number of $A$-type smoothings of $s$ minus the number of $B$-type smoothings of $s$. Observe that $s$ has no crossings and has the same number of segment constituents as $D$, though in general, there is no natural correspondence between the segment constituents of $D$ and $s$.

For every state $s$, we define auxiliary polynomials $G_s$, $L_s$, and $E_s$ to record the data of the state's underlying graph, loop constituents, and segment constituents, respectively. They are as follows.

Introduce variables $\lambda_{\{P,Q\}}$ for every unordered pair $\{P, Q\}$ of not necessarily distinct poles $P,Q \in P(\diag)$ and let $s(P,Q)$ denote the number of edges between $P$ and $Q$ in the underlying graph of $s$. Define the polynomial
\[G_s = \prod_{\{P,Q\}}\lambda_{\{P, Q\}}^{s(P,Q)}.\]

Let $2^{P(\diag)}/\sim$ denote the set of unordered bipartitions of the poles of $\diag$. (The equivalence relation $\sim$ identifies each subset of $P(\diag)$ with its complement.) Each loop constituent of $s$ separates the sphere into two regions and determines an element of $2^{P(\diag)}/\sim$. For each $U \in 2^{P(\diag)}/\sim$, let $n_U(s)$ denote the number of loop constituents of $s$ corresponding to $U$. Let $n_0(s) := n_{[\emptyset]}(s)$ denote the number of loops that create a region with no poles. We call such loops \emph{nullhomotopic}. We define the polynomial $L_s$ in variables $\{x_U \mid U \in 2^{P(\diag)}/\sim, U \neq [\emptyset]\}$ by
\[L_s = \prod_{\substack{U \in 2^{P(\diag)}/\sim\\U \neq [\emptyset]}}x_U^{n_U(s)}.\]

Fix a set of shortcuts $\{\alpha_{P,Q} \mid P, Q \in P(\diag)\}$ between every pair of poles of $\diag$, where $\alpha_{P,Q}$ is oriented from $P$ to $Q$. Suppose $e, f$ are segment constituents, each belonging to either $\diag$ or one of its states. Given orientations on $e$ and $f$, define
\[\mu_{e,f}(P,Q) = e \cdot \alpha_{P,Q} - f \cdot \alpha_{P,Q}.\]

We say $e$ and $f$ are \emph{aligned} and write $e \parallel f$ if their endpoint poles coincide and they are identical (as unoriented constituents) in a neighborhood of each of the endpoint poles.

Suppose $e \in E(\diag)$ is oriented from pole $e_0$ to pole $e_1$. For a given state $s$, there is at most one $e_s \in E(s)$ aligned with $e$. If such an $e_s$ exists, we give it the orientation that agrees with $e$ near $e_0$ and $e_1$. Define the polynomial $E_s$ in variables $\{y_{e, P} \mid e \in E(\diag), P \in P(\diag)\}$ by
\[E_s = \prod_{\substack{e \parallel e_s\\p \in P(\diag)}}y_{e, P}^{\mu_{e_s, e}(e_0, P)},\]
where the product is taken over all $e \in E(\diag), P \in P(\diag)$ such that $s$ has a segment constituent $e_s \in E(s)$ with $e \parallel e_s$.

Finally, define the \emph{generalized bracket polynomial} $\llangle \diag \rrangle$ as the Laurent polynomial in variables $A$, $\{\lambda_{\{P,Q\}}\}$, $\{x_U\}$, and $\{y_{e, P}\}$ given by
\begin{equation}\label{bracketdef}
    \llangle \diag \rrangle = \sum_{s} A^{\sigma(s)}(-A^2-A^{-2})^{n_0(s)} \cdot G_s \cdot L_s \cdot E_s.
\end{equation}
We show in Lemma \ref{path-independence-two} that $E_s$ is well-defined, independent of the choice of shortcuts $\{\alpha_{P,Q}\}$. The lemma should be thought of as an analog of Corollary \ref{path-independence} for closed (but not generic) curves. We argue similarly to the proof given for Lemma 8.1 in \cite{Turaev1}. (See also the \emph{shortcut moves} of \cite{linov}.)
\begin{lemma}\label{path-independence-two}
    Suppose $e, f$ are segment constituents, each belonging to either $\diag$ or one of its states, with $e \parallel f$ and orientations that agree near their endpoints. Let $P, Q \in P(\diag)$. Then $\mu_{e,f}(P,Q)$ does not depend on the choice of shortcut $\alpha_{P,Q}$.
\end{lemma}
\begin{proof}
     Observe that any shortcut between $P$ and $Q$ is obtainable from $\alpha_{P,Q}$ by a sequence of local transformations of the following types. 
    \begin{enumerate}[(i)]
        \item Pulling $\alpha_{P,Q}$ across a strand of $\diag$ (c.f. Type II Reidemeister),
        \item Pulling $\alpha_{P,Q}$ across a crossing of $\diag$ (c.f. Type III Reidemeister),
        \item Pulling $\alpha_{P,Q}$ across a pole other than $P$ or $Q$ (c.f. pole slide move),
        \item Pulling $\alpha_{P,Q}$ across a strand of $\diag$, near $P$ or $Q$ (c.f. pole twist move).
    \end{enumerate}
    Observe that (i) and (ii) each preserve $e \cdot \alpha_{P,Q}$ and $f \cdot \alpha_{P,Q}$, while (iii) and (iv) each preserve $e \cdot \alpha_{P,Q} - f \cdot \alpha_{P,Q}$ by the alignment condition.
\end{proof}
\begin{corollary}\label{muzero}
    With $e,f$ as in Lemma \ref{path-independence-two}, we have $\mu_{e,f}(P,P) = 0$.
\end{corollary}
The \emph{normalized generalized bracket polynomial} $\llangle \diag \rrangle_*$ is defined by
\[\llangle \diag \rrangle_* = (-A^3)^{-\on{wr}(\diag)} \llangle \diag \rrangle.\]
Here, the \emph{writhe} $\on{wr}(\diag)$ of the oriented generalized knotoid diagram $\diag$ is defined by the sum of $\on{sgn}(c)$ over all crossings $c$ in $\diag$.
\begin{proposition}
    The polynomial $\llangle \diag \rrangle_*$ is a generalized knotoid invariant.
\end{proposition}
\begin{proof}
    Observe that the generalized bracket polynomial satisfies the usual disjoint union and skein relations
    \[\llangle \diag \cup \bigcirc \rrangle = (-A^2-A^{-2})\llangle \diag \rrangle \quad \text{and} \quad \llangle \diag \rrangle = A\llangle \diag_+ \rrangle + A^{-1}\llangle \diag_- \rrangle,\]
    where $\diag_+$ and $\diag_-$ denote the diagrams obtained by performing an $A$-smoothing and a $B$-smoothing, respectively, at a particular crossing of $\diag$. As with the traditional bracket polynomial for knots, this implies that the normalized polynomial $\llangle \diag \rrangle_*$ is preserved by Reidemeister moves.
\end{proof}
In analogy to Proposition \ref{heightbound} for the base-pointed index polynomial, Propsition \ref{bracketheightbound} gives lower bounds on heights between poles in terms of the generalized bracket polynomial. Before stating the bounds, we record a lemma.

\begin{figure}[htbp]
    \centering
    \includegraphics{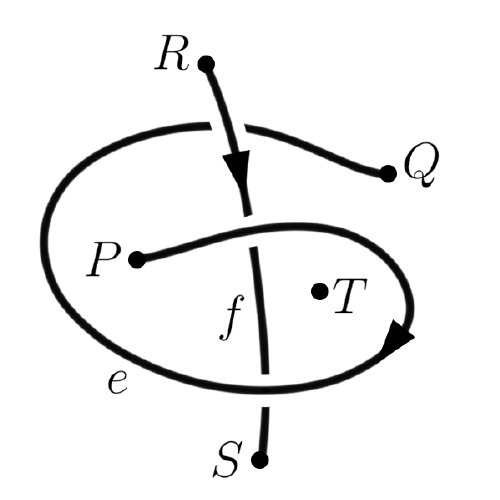}
    \caption{An oriented generalized knotoid with five poles $P,Q,R,S,T$ and two constituents $e,f$.}
    \label{bracketpolyexample}
\end{figure}

\begin{lemma}\label{triangle}
    Let $e$,$f$ be as in Lemma \ref{path-independence-two} and let $P,Q,R \in P(\diag)$. Then
    \[\mu_{e,f}(P,Q) + \mu_{e,f}(Q,R) = \mu_{e,f}(P,R).\]
\end{lemma}
\begin{proof}
    If $P = Q$ or $Q = R$, the result follows from Corollary \ref{muzero}, so assume this is not the case. Without loss of generality, we choose $\alpha_{P,Q}$ and $\alpha_{Q,R}$ with no intersections other than poles; by Lemma \ref{path-independence-two}, this does not change any of the values in question. Then we obtain a shortcut $\tilde \alpha_{P,R}$ by concatenating $\alpha_{P,Q}$ with $\alpha_{Q,R}$ and perturbing the resulting curve in a neighborhood of $Q$ so that, locally, it neither passes through $Q$ nor intersects $e$ or $f$. (This is possible even if $e$ and $f$ have an endpoint at $Q$ since $e \parallel f$.) Observe that $e \cdot \alpha_{P,Q} + e \cdot \alpha_{Q, R} = e \cdot \tilde \alpha_{P,R}$ and $f \cdot \alpha_{P,Q} + f \cdot \alpha_{Q, R} = f \cdot \tilde \alpha_{P,R}$, whence
    \begin{align*}
        \mu_{e,f}(P,Q) + \mu_{e,f}(Q,R) &= (e \cdot \alpha_{P,Q} - f \cdot \alpha_{P,Q}) + (e \cdot \alpha_{Q,R} - f \cdot \alpha_{Q,R})\\
        &= e \cdot \tilde \alpha_{P,Q} - f \cdot \tilde \alpha_{P,Q}\\
        &= \mu_{e,f}(P,R),
    \end{align*}
    where the last equality follows from Lemma \ref{path-independence-two}.
\end{proof}
Following the notation of \cite{Turaev1}, we define the \emph{$x_i$-span}, denoted $\on{spn}_{x_i}(F)$, of a Laurent polynomial $F \in \mathbb Z[x_1^{\pm 1}, \dots, x_n^{\pm 1}]$ as the non-negative difference between the largest and smallest exponents of $x_i$ among monomials with nonzero coefficients in $F$ (by convention, $\on{spn}_{x_i}(0) = -\infty$). Given a generalized knotoid $\kappa$ and poles $P, Q \in P(\kappa)$, we let $\{P \mid Q\}$ denote the set of variables $x_U$ with $U \in 2^{P(\kappa)}/\sim$ such that $P$ and $Q$ are separated by $U$.
\begin{proposition}\label{bracketheightbound}
    Let $\kappa$ be a pole-labeled, constituent-labeled, oriented generalized knotoid and let $P, Q \in P(\kappa)$ be poles. Then
    \begin{enumerate}[(a)]
        \item $h_\kappa(P, Q) \ge \deg_{\{P \mid Q\}}(\llangle \kappa \rrangle_*)$,
        \item $h_\kappa(P, Q) \ge \frac 12\max_{e \in E(\kappa)} \deg_{y_{e, P}}\left(\llangle \kappa\rrangle_* \Big|_{y_{e, Q} = y_{e, P}^{-1}}\right)$,
        \item $h_\kappa(P, Q) \ge \frac 12\max_{e \in E(\kappa)} \on{spn}_{y_{e, P}}\left(\llangle \kappa\rrangle_* \Big|_{y_{e, Q} = y_{e, P}^{-1}}\right)$.
    \end{enumerate}
\end{proposition}
\begin{proof}
    Let $\diag$ be any diagram of $\kappa$ and let $\gamma$ be any shortcut oriented from $P$ to $Q$. Let $d_a$ denote the right side of the inequality in (a). Then there is a state $s$ of $\diag$ with $d_a = \sum_{U \in \{P\mid Q\}} n_U(s)$. Since $\gamma$ has at least one intersection with each loop constituent of $s$ that separates $P$ and $Q$, part (a) follows.
    
    We now prove (b) and (c). First, suppose $e,f$ are oriented segment constituents, each belonging to either $\diag$ or a state of $\diag$. Each point of intersection of $\gamma$ with $e \cup f$ away from poles is counted by $e \cdot \gamma - f \cdot \gamma$ with multiplicity in $\{-2, -1, 0, 1, 2\}$. It follows that the number of such points, denoted $\#(\gamma \cap (e \cup f))$, satisfies
    \begin{equation}\label{bound}
        \#(\gamma \cap (e \cup f)) \ge \frac 12 |e \cdot \gamma - f \cdot \gamma|.
    \end{equation}
    Let $d_b$ denote the right side of the inequality in (b). Suppose $d_b > 0$. Then there is a segment constituent $e \in E(\diag)$, a state $s$, and a segment constituent $e_s \in E(s)$ such that $e \parallel e_s$ and
    \begin{align*}
        d_b &= \frac 12 |\mu_{e_s, e}(e_0, P) - \mu_{e_s, e}(e_0, Q)|\\
        &= \frac 12 |\mu_{e_s, e}(P, Q)|\\
        &= \frac 12 |e_s \cdot \gamma - e \cdot \gamma|,
    \end{align*}
    where the second equality uses Lemma \ref{triangle} and the third uses Lemma \ref{path-independence-two}. It follows from Inequality~\eqref{bound} that $\#(\gamma \cap (e_s \cup e)) \ge d_b$, and part (b) follows.
    
    Finally, let $d_c$ denote the right side of the inequality in (c). Assume $d_c > d_b$. Then there is a segment constituent $e \in E(\diag)$, two states $s_1, s_2$, and two segment constituents $f_1 \in E(s_1)$, $f_2 \in E(s_2)$ such that $e \parallel f_1 \parallel f_2$ and
    \begin{align*}
        d_c &= \frac 12 |(\mu_{f_1, e}(e_0, P) - \mu_{f_1, e}(e_0, Q)) - (\mu_{f_2, e}(e_0, P) - \mu_{f_2, e}(e_0, Q))|\\
        &= \frac 12 |-\mu_{f_1, e}(P, Q) + \mu_{f_2, e}(P,Q)|\\
        &= \frac 12 |\mu_{f_2, f_1}(P,Q)|\\
        &= \frac 12 |f_2 \cdot \gamma - f_1 \cdot \gamma|,
    \end{align*}
    where the second equality uses Lemma \ref{triangle} and the fourth uses Lemma \ref{path-independence-two}. It follows from Inequality~\eqref{bound} that $\#(\gamma \cap (f_1 \cup f_2)) \ge d_c$, and part (c) follows.
\end{proof}
\begin{example}
    Let $\kappa$ denote the generalized knotoid represented by the diagram in Figure \ref{bracketpolyexample}. Computations yield
    \begin{align*}
        \llangle \kappa \rrangle_* = &-\lambda_{QR}\lambda_{PS}(1 + A^{-2} + A^{-2}x_{[\{T\}]} + A^{-4}) - \lambda_{PR}\lambda_{QS}A^{-4}\\
        &-\lambda_{PQ}\lambda_{RS}y_{e,Q}y_{e,R}y_{e,S}y_{f,P}^{-1}(A^{-2}y_{f,T} + A^{-4}x_{[\{T\}]} + A^{-6}y_{e,T}^{-1}).
    \end{align*}
    For brevity, we have omitted the set notation in the subscripts of the $\lambda$ variables. Proposition \ref{bracketheightbound}(a) yields the height bounds $h_{\kappa}(T, P)$, $h_{\kappa}(T, Q)$, $h_{\kappa}(T, R)$, $h_{\kappa}(T, S) \ge 1$. Meanwhile, Proposition \ref{bracketheightbound}(b) and \ref{bracketheightbound}(c) yield the height bounds $h_\kappa(P, Q)$, $h_\kappa(P, R)$, $h_\kappa(P, S) \ge 1/2$ and $h_{\kappa}(T, P)$, $h_{\kappa}(T, Q)$, $h_{\kappa}(T, R)$, $h_{\kappa}(T, S) \ge 1$. In this case, the bounds give enough information for a complete description of the height spectrum of $\kappa$.
\end{example}
We remark on the case of generalized knotoids that are not necessarily oriented or fully labeled. Consider an oriented, pole-labeled, constituent-labeled diagram $\diag$. As we remove orientation or labeling data from $\diag$, we alter the polynomial $\llangle \diag \rrangle_*$ accordingly, as follows.
\begin{enumerate}[(i)]
    \item If constituent labels are removed, we re-index the variables $\{y_{e,P}\}$ by replacing each subscript component $e \in E(\diag)$ with the ordered pair $(e_0, e_1)$ corresponding to its endpoint poles.
    \item If pole labels are removed, we remove $G_s$ and $E_s$ from the definition given in Equation~\eqref{bracketdef}, and we re-index the variables $\{x_U\}$ as follows. For each subscript $U$, write $U = [S]$ for a set $S \subset P(\diag)$, then replace $U$ with the integer $\min(|S|, |P(\diag)| - |S|)$.
    \item If orientation is disregarded, we remove $E_s$ from the definition given in Equation~\eqref{bracketdef} and redefine the normalized bracket polynomial by $\llangle \diag \rrangle_* = (-A^3)^{-\sum_{e \in C(\diag)}\on{wr}(e)}\llangle \diag \rrangle$, where $\on{wr}(e)$ denotes the sum of the signs of the self-crossings of an arbitrarily-oriented constituent $e$. (In particular, $\on{wr}(e)$ is independent of orientation.)
\end{enumerate}
Given an oriented, unlabeled diagram $\diag$ (for example), we obtain the bracket polynomial $\llangle \diag \rrangle_*$ by labeling $\diag$ arbitrarily and then applying (i) and (ii). (Note that (ii) and (iii) subsume (i).) It is straightforward to verify that, in general, this method yields a bracket polynomial invariant for each class of generalized knotoids with specified labeling and orientation data.
\subsection{Recovering existing bracket polynomials}
We survey some existing bracket polynomial constructions for classes of objects subsumed by generalized knotoids and show that they are recovered by the generalized bracket polynomial. 

In the case of a classical link diagram $L$, there are no poles, so all loop constituents of any state $s$ are nullhomotopic and $n_0(s)$ is simply the number of constituents of $s$. Moreover, $G_s$, $L_s$, and $E_s$ are constant, so $\llangle L \rrangle_*$ recovers the traditional bracket polynomial $\langle L \rangle$ of the link diagram.

In \cite{multi-linkoid}, a bracket polynomial for oriented, ordered (i.e. pole-labeled) multi-linkoids $L$ is defined by
\[\langle L \rangle_\bullet(A, \{\lambda_{ij}\}) = (-A^3)^{-\on{wr}(L)}\sum_{s} A^{\sigma(s)}(-A^2-A^{-2})^{\|s\|}\prod_{\Lambda} \lambda_{ij}.\]
Here, the poles of the multi-linkoid are labeled $1, 2, \dots, 2n$, the symbol $\Lambda$ indicates that the product is taken over all $i < j$ such that poles $i$ and $j$ are connected by a segment constituent in $s$, and $\|s\|$ denotes the number of loop constituents of $s$. It is clear that $\langle L \rangle_\bullet$ is recovered from $\llangle L \rrangle_*$ by identifying $G_s$ with $\prod_{\Lambda}\lambda_{ij}$ and making the substitutions $x_U = -A^2 - A^{-2}$ for all $U \in (2^{P(\diag)}/\sim)\setminus [\emptyset]$ and $y_{e, P} = 1$ for all $e \in E(L), P \in P(L)$.

In \cite{Turaev1}, Turaev defines a two-variable bracket polynomial for oriented spherical knotoids. Let $K$ be a spherical knotoid oriented from $L$ to $H$, and let $\alpha = \alpha_{L,H}$ be a shortcut oriented from $L$ to $H$. (We also use $K$ to denote the segment constituent of the knotoid.) Turaev's polynomial is given by
\[\llangle K \rrangle_\circ(A, u) = (-A^3)^{-\on{wr}(K)}\sum_{s}A^{\sigma(s)}(-A^2 - A^{-2})^{\|s\|}u^{k_s \cdot \alpha - K \cdot \alpha},\]
where $k_s$ denotes the unique segment constituent of the state $s$, oriented from $L$ to $H$. Note that for any state $s$, the segment constituent $k_s$ aligns with $K$, and all loop constituents of $s$ are nullhomotopic. In terms of the generalized bracket polynomial, the polynomials $G_s = \lambda_{LH}$ and $L_s = 1$ are constant, and $n_0(s) = \|s\|$. It follows that $\llangle K \rrangle_\circ$ is recovered from $\llangle K \rrangle_*$ by the substitutions $\lambda_{LH} = 1$  and $y_{K, H} = u$. It also follows that the height bound in Proposition \ref{bracketheightbound}(c) is equivalent to the complexity bound given in Equation 8.3.1 of \cite{Turaev1}.

In \cite{Turaev1}, Turaev extends $\llangle K \rrangle_\circ$ to a bracket polynomial for oriented planar knotoids given by
\[[K]_\circ(A, B, u) = (-A^3)^{-\on{wr}(K)} \sum_s A^{\sigma(s)}(-A^2 - A^{-2})^{p(s)}B^{q(s)}u^{k_s \cdot \alpha - K\cdot \alpha},\]
where $p(s)$ denotes the number of loop constituents of $s$ that do not enclose the the segment constituent $k_s$ and $q(s)$ denotes the number of loop constituents of $s$ that enclose $k_s$. Represent $K$ as a generalized knotoid with poles $L,H, \infty$ and constituent $K$ oriented from $L$ to $H$, equipped with shortcuts $\{\alpha_{P,Q} \mid P,Q \in \{L, H, \infty\}\}$ where $\alpha_{L,H} = \alpha$. Then $p(s) = n_0(s)$ is simply the number of nullhomotopic loop constituents of $s$, and $q(s) = n_{[\{\infty\}]}(s)$ is the number of loop constituents of $S$ that determine the bipartition $\{\{\infty\}, \{L, H\}\}$. Thus $[K]_\circ$ is recoverable from $\llangle K \rrangle_*$ by the substitutions $\lambda_{LH} = 1$, $x_{[\{\infty\}]} = B$, $y_{K, \infty} = 1$, and $y_{K, H} = u$.

In \cite{kutluay}, Kutluay refines $[K]_\circ$ by replacing the term $u^{k_s \cdot \alpha - K \cdot \alpha}$ with the term $\ell^{w_{\gamma_s}(L) - w_\gamma(L)}h^{w_{\gamma_s}(H) - w_\gamma(H)}$, obtaining a polynomial $[K]_{\bullet}(A,B,\ell, h)$. Here, $w_\beta \colon \mathbb R^2 \to \frac 12 \mathbb Z$ is the \emph{winding potential function} of an oriented generic closed curve $\beta$ that maps each point $x \in \mathbb R^2$ to the winding number of $\beta$ around $x$. (For $x \in \beta$, the winding number is taken to be the average of the winding numbers of the regions adjacent to $x$.) The closed curve $\gamma$ (resp. $\gamma_s$) is defined by $K \cup \alpha^r$, (resp. $k_s \cup \alpha^r$). Kutluay observes that
\[k_s \cdot \alpha - K \cdot \alpha = [w_{\gamma_s}(L) - w_\gamma(L)] - [w_{\gamma_s}(H) - w_\gamma(H)],\]
so Turaev's polynomial $[K]_\circ$ is recoverable from Kutluay's polynomial $[K]_\bullet$ via the substitutions $\ell = u$ and $h = u^{-1}$.

It is well-known that for a generic closed curve $\beta$ on the plane and a point $x \not \in \beta$, the winding number $w_\beta(x)$ is given by $\beta \cdot \alpha_{x, \infty}$, where $\alpha_{x,\infty}$ is any path from $x$ to the unbounded region of the plane that intersects $\beta$ transversely and away from self-intersections of $\beta$. Using this fact, it is straightforward to verify that
\begin{align*}
    w_{\gamma_s}(L) - w_\gamma(L) &= \gamma_s \cdot \alpha_{L, \infty} - \gamma \cdot \alpha_{L, \infty}\\
    &= \mu_{k_s, K}(L, \infty),
\end{align*}
and similarly
\begin{align*}
    w_{\gamma_s}(H) - w_\gamma(H) &= \mu_{k_s, K}(H, \infty)\\
    &= -\mu_{k_s, K}(L, H) + \mu_{k_s, K}(L, \infty).
\end{align*}
It follows that $[K]_\bullet$ is recoverable from $\llangle K \rrangle_*$ by the substituions $\lambda_{LH} = 1$, $x_{[\{\infty\}]} = B$, $y_{K, \infty} = \ell h$, and $y_{K, H} = h^{-1}$.

\begin{comment}
\section{An Alexander Polynomial for Staked Multiknotoids}\label{bracket}
\subfile{./sections/alexpoly.tex}
\end{comment}

\section{Knotoidal Graphs}\label{knotoidal}
\subsection{Knotoidal graphs}
We now define an extension of generalized knotoids. Let $\Sigma$ denote a closed orientable surface, and let $G$ be a finite graph. As with generalized knotoids, we do not require $G$ to be connected or simple, and $G$ may have valency-zero vertices. Let $\Tilde{G}$ denote the disjoint union of $G$ with a finite collection of circles. The edges and circles of $\Tilde{G}$ are still called constituents.

\begin{figure}[htbp]
    \centering
    \includegraphics[]{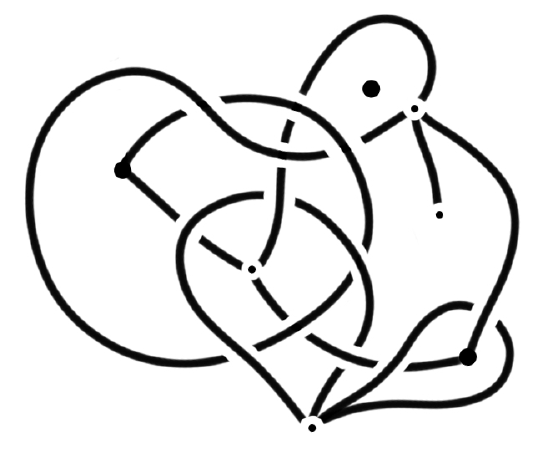}
    \caption{A knotoidal graph with three poles and four spatial vertices.}
    \label{knotoidalgraphexample}
\end{figure}

A \textit{knotoidal graph diagram} on $\Sigma$ is a pair $(\mathcal{D}, V(\mathcal{D}))$, where $\mathcal{D}$ is a generic immersion of $\Tilde{G}$ in $\Sigma$ whose only singularities are transverse double points, called crossings, with over/undercrossing data, and where $V(\mathcal{D})$ is a chosen subset of the images of the vertices of $G$, each of which has valency at least one. We call the elements of $V(\mathcal{D})$ the \textit{spatial vertices} of $\mathcal{D}$. The images of the remaining vertices of $\Tilde{G}$ are called \textit{poles} of $\mathcal{D}$ and the set of poles is denoted by $P(\diag)$. Figure \ref{knotoidalgraphexample} shows an example.

The rest of the terminology for generalized knotoids carries over. we call the graph $G$ (resp. $\Tilde{G}$) the \textit{underlying graph} (resp. \textit{underlying looped graph}) of $\mathcal{D}$. The \textit{valency} of a pole or spatial vertex is the valency of the corresponding vertex in the underlying graph $G$. Let $E(\diag)$ denote the set of images of the edges of $\Tilde{G}$, called \textit{segment constituents} of $\mathcal{D}$, and let $L(\mathcal{D})$ denote the images of the circles of $\Tilde{G}$, called \textit{loop constituents of $\diag$}. Let $C(\diag) := E(\diag) \cup L(\diag)$ denote the set of \textit{constituents of $\diag$}. 

\begin{figure}[htbp]
    \centering
    \includegraphics[width=12cm]{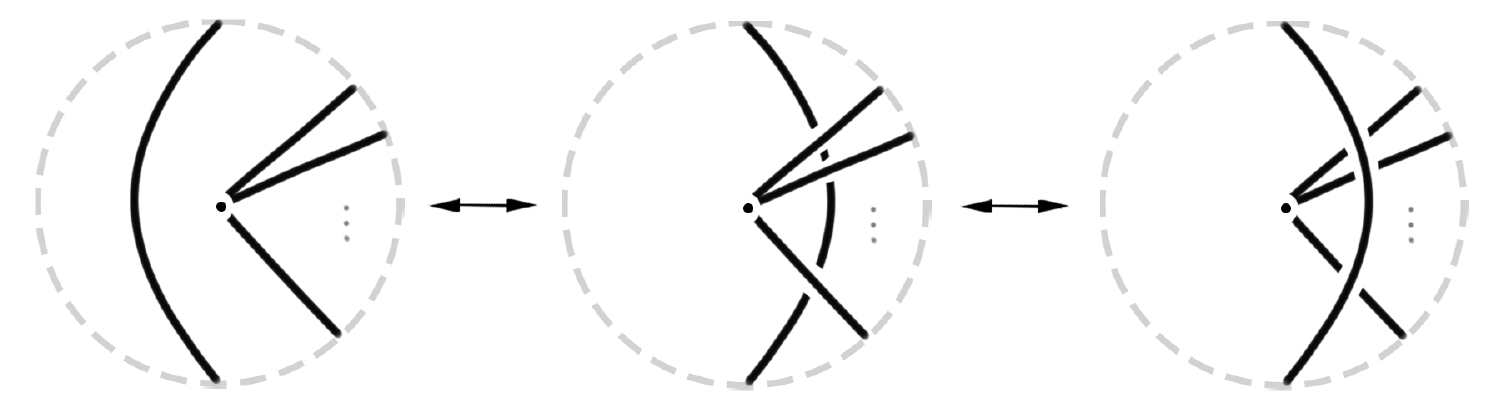}
    
    \includegraphics[width=12cm]{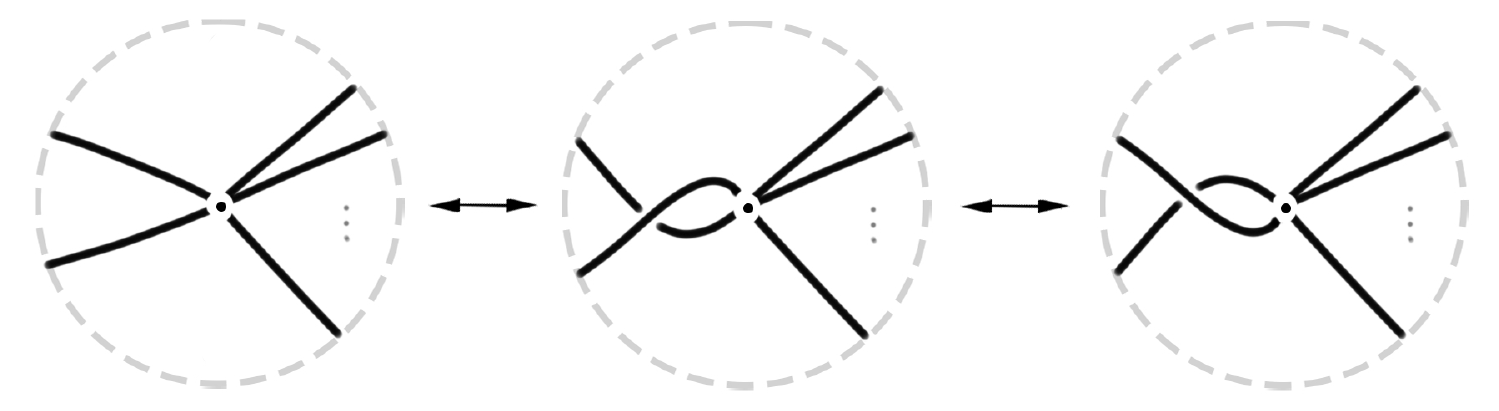}
    \caption{The vertex slide move (top) and the vertex twist move (bottom).}
    \label{genRMoves}
\end{figure}

We consider knotoidal graph diagrams in $\Sigma$ up to ambient isotopy and \textit{generalized Reidemeister moves}: the three standard Reidemeister moves away from poles and spatial vertices, along with the \emph{vertex slide move} and the  \emph{vertex twist move} near spatial vertices. See Figure \ref{genRMoves}. The forbidden pole moves of Figure \ref{forbidden} remain in effect. A \textit{knotoidal graph} $KG$ is an equivalence class of knotoidal graph diagrams. We may speak of knotoidal graphs possibly with labels assigned to poles or spatial vertices, and labels or orientations assigned to constituents.

With the vertex twist move, the theory of knotoidal graphs is akin to the theory of spatial graphs with pliable vertices (as opposed to rigid vertices, where the twist move is not allowed). One could require knotoidal graphs to have rigid spatial vertices, which would yield a different theory that we do not consider in this paper.

As we remarked for generalized knotoids, an equivalent theory arises by considering knotoidal graphs in compact surfaces with boundary while disallowing valency-zero poles. 

\begin{example}
Generalized knotoids are knotoidal graphs with no spatial vertices.
\end{example}

\begin{example}
Spatial graphs are knotoidal graphs with no poles. 
\end{example}

\begin{example}
A \emph{graphoid} (resp. \emph{multi-graphoid}) as defined in \cite{bonded} is a knotoidal graph with exactly two poles (resp. $2n$ poles for some $n \ge 1$), each of valency one. A \emph{bonded knotoid}, used to model bonded proteins \cite{bondedknotoids},  is particular type of graphoid obtained by adding edges to a knotoid with endpoints at distinct points of the knotoid away from its endpoints. So all resulting spatial vertices have valency three.
\end{example}
\subsection{Rail diagrams}\label{railsection}
\begin{figure}[hbtp]
\centering
\includegraphics[scale=0.8]{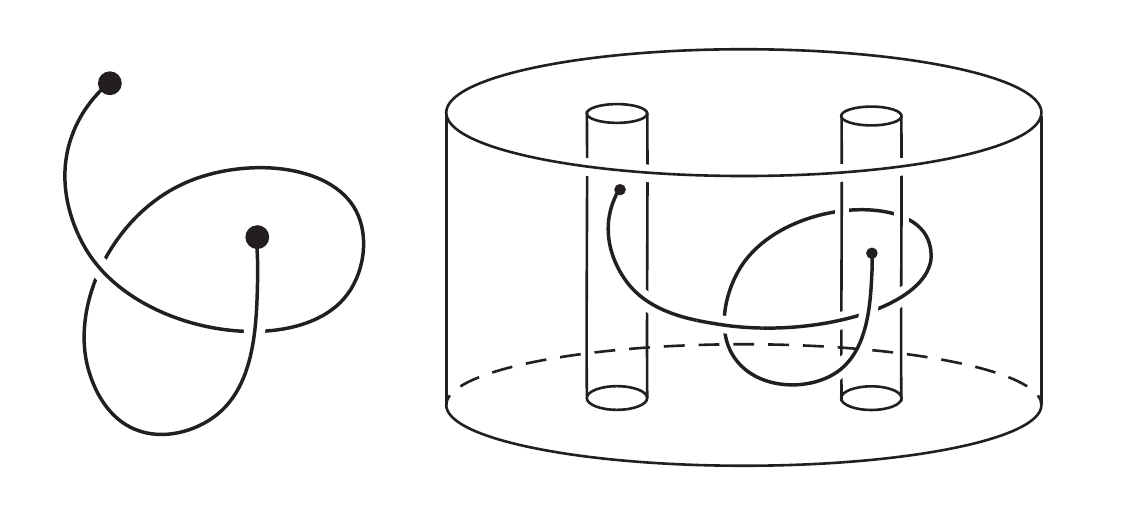}
\caption{A planar rail diagram for a knotoid.}
\label{raildiagram}
\end{figure}

A topological perspective on knotoids is introduced in \cite{rails} by means of \emph{rail diagrams}. For planar knotoids, consider the thickened disk $D^2 \times I$ together with two distinct lines  $\ell_1 = \{x\} \times I, \ell_2 = \{y\} \times I$, called \textit{rails}. A \emph{planar rail diagram} is a proper embedding of an arc into $(D^2 \times I) \setminus \mathring{N}(\ell_1 \cup \ell_2)$ such that one endpoint is embedded in $\partial N(\ell_1)$ and one is embedded in $\partial N(\ell_2)$, as shown in Figure \ref{raildiagram}. We consider planar rail diagrams up to isotopies of the embedded arc in $(D^2 \times I) \setminus \mathring{N}(\ell_1 \cup \ell_2)$. Intuitively, the removal of the rails has the effect of topologically enforcing the forbidden move for knotoid diagrams. It is shown in \cite{rails} that the theory of planar rail diagrams is equivalent to the theory of planar knotoids, and if the disk $D^2$ is replaced by an arbitrary surface $\Sigma$, the arguments in \cite{rails} generalize readily to show that the theory of $\Sigma$-rail diagrams is equivalent to the theory of knotoids on $\Sigma$.

We extend $\Sigma$-rail diagrams to knotoidal graphs. For a given knotoidal graph $KG$ in $\Sigma$ with diagram $\mathcal{D}$, take the thickened surface $\Sigma \times I$ and $n$ distinct rails $\ell_i = \{x_i\} \times I$, $i = 1, \dotsc, n$, where $n = |P(\mathcal{D})|$, each rail corresponding to a pole $p_i$. Let $\Tilde{G}'$ be the graph obtained from the underlying looped graph $\Tilde{G}$ of $\mathcal{D}$ by removing a disk neighborhood of each pole. Then for each pole $p_i$ there are $v(p_i)$ endpoints in $\Tilde{G}'$ created by removing the disk neighborhood. Now properly embed $\Tilde{G}'$ in $(\Sigma \times I) \setminus \mathring{N}(\bigcup_{i = 1}^n \ell_i)$ so that the projection of the embedding to $\Sigma \times \{0\}$ recovers $\mathcal{D}$ away from the poles, and so the $v(p_i)$ endpoints of $\Tilde{G}'$ are embedded on the boundary $\partial N(\ell_i)$ of the corresponding rail neighborhood. 

Again, we consider these diagrams up to isotopies of the embedded looped graph in $(\Sigma \times I) \setminus \mathring{N}(\bigcup_{i = 1}^n \ell_i)$, with one major difference. We may identify the cylindrical boundary portion of each $\partial N(\ell_i)$ with $S^1 \times I$. If an endpoint $v$ of $\Tilde{G}'$ is embedded as the point $(x, t)$ in $\partial N(\ell_i)$, then we require the isotopies to keep $v$ in the vertical line segment $\{x\} \times I$. That is, we only allow endpoints of the graph $\Tilde{G}'$ embedded in rail neighborhoods to slide vertically up and down. This prevents arbitrary twisting near poles, which we disallow in the diagrammatic theory of knotoidal graphs. 

The theory of $\Sigma$-rail diagrams of knotoidal graphs is equivalent to the theory of knotoidal graphs in $\Sigma$. Consequently, we may derive invariants of knotoidal graphs using topological methods, as we do in Section \ref{hyperbolicity}.

Although we do not address it here, one can also allow virtual crossings in a spherical or planar knotoidal graph. This would generate classical knotoidal graphs in surfaces as we have described, however, in this theory, different projections can generate surfaces of different genus and equivalence between diagrams must then allow adding and removing handles for the corresponding surfaces. See \cite{graphoids} for this theory as applied to graphoids.

\section{Hyperbolicity for Knotoidal Graphs}\label{hyperbolicity}
We review the notion of hyperbolicity for knotoids, studied in \cite{hypknotoids}. Two maps $\phi_{S^2}^D$ and $\phi_{S^2}^G$ are defined in \cite{hypknotoids}, each sending a spherical knotoid to a knot in the thickened torus $T \times (0, 1)$. Two maps $\phi_{\mathbb{R}^2}^D$ and $\phi_{\mathbb{R}^2}^G$ are also defined in \cite{hypknotoids}, each sending a planar knotoid to a knot in the genus three handlebody $H_3$ or to a knot in the genus two handlebody $H_2$, respectively. 

The map $\phi_{S^2}^D$ is termed the \emph{spherical reflected doubling map} and is defined as follows.

%Throughout our paper, we will find it useful to perform an operation on a 3-manifold $M$ with boundary and a specified subsurface $F$ of the boundary that involves taking a reflected copy of our manifold $M'$ with corresponding specified subsurface $F'$ and identifying each point in $F$ with the corresponding point in $F'$. We define it below as the \emph{reflected double}. {\color{teal}We don't actually use this until the very last section, right? It might be better if we move this paragraph to the beginning of that section} {\color{magenta} probably needs a diagram too}

\begin{definition} \label{spherical_doubling_map}
Consider a spherical rail diagram for a knotoid $k$, namely a thickened cylinder $M = S^2 \times (0, 1) \setminus \mathring{N}(\ell_1 \cup \ell_2)$ with an embedded arc, where each $\ell_i = \{x_i\} \times (0, 1)$ is a rail. Let $C_i = \partial N(\ell_i)$ denote the boundaries of the regular neighborhoods of the rails, which are homeomorphic to annuli and are each punctured once by the knotoid (we may take the neighborhoods to be sufficiently small so that this is true). Label the punctures by $z_i$. Take a reflected copy $M^R$ of $M$ which contains reflected copies $C_i^R$ of the rail boundaries and copies $z_i'$ of the punctures. Then we glue together $M$ and $M^R$ by gluing $C_1$ to $C_1^R$ and $C_2$ to $C_2^R$ via a reflection. This yields a well-defined knot $K$ in the thickened torus $T \times (0, 1)$. We then set $\phi_{S^2}^D$ to $K \in \mathcal{K}(T^2 \times (0, 1))$.
\end{definition}

We say that $k$ is hyperbolic if $K$ is hyperbolic in $T \times (0, 1)$, that is, $(T \times (0, 1)) \setminus N(K)$ admits a complete hyperbolic metric. We define the volume of $k$, denoted $\vol_{S^2}(k)$,  to be half the volume of $(T \times (0, 1)) \setminus N(K)$. 

\subsection{Hyperbolicity for Knotoidal Graphs}
We extend the domain of the map $\phi_{S^2}^D$ to knotoidal graphs.

\begin{definition} \label{defnKGhyp}
Consider a spherical rail diagram for a knotoidal graph $KG$ in $S^2$. Denote by $\ell_1, \dotsc, \ell_m$ the rails corresponding to nonzero-valency poles $p_1, \dotsc, p_m$, and denote by $\ell_{m+1}, \dotsc, \ell_{n}$ the rails corresponding to valency-zero poles. Let $C_i$ be the boundary of $N(\ell_i)$ for $i = 1, \dotsc, m$; we may choose the regular neighborhoods to be sufficiently small such that each $C_i$ is punctured exactly $v(p_i)$ times, where $v(p_i)$ is the valency of the pole $p_i$. Let $M$ be the manifold $(S^2 \times I) \setminus N(\bigcup_{i = 1}^n \ell_i)$, and take a reflected copy $M^R$ which contains reflected copies $C_i^R$ of the rail boundaries. Then glue together $M$ and $M^R$ by gluing $C_i$ to $C_i^R$ via reflection for $i = 1, \dotsc, m$, so that the punctures in the $C_i$ and $C_i^R$ line up. Note that we do not glue the boundaries corresponding to valency-zero poles, and if there are no poles of nonzero valency, then we do not take a reflected copy, leaving us with just $M$. 

This process yields a spatial graph $\mathcal{G}$ in a manifold $Y$. If there are no valency-zero poles, then $Y$ is homeomorphic to a thickened closed orientable surface of genus $m - 1$. Otherwise, $Y$ is homeomorphic to a handlebody of genus $m + (m - 1) + 2(n - m - 1) = 2n - 3$ when $n - m > 1$: each valency-zero pole after the first one increases the genus by two. If there is exactly one valency-zero pole, then $Y$ is homeomorphic to a handlebody of genus $m + (m - 1) = 2m - 1$.  Hence, we obtain an extension of the \textit{spherical reflected doubling map} $\phi_{S^2}^D$ to knotoidal graphs. It associates to every knotoidal graph in $S^2$ a spatial graph in a 3-manifold $Y$, and we write $\phi_{S^2}^D(KG) = (\mathcal{G}, Y)$. 
\end{definition}

When restricted to spherical knotoids (knotoidal graphs with two poles and a single segment constituent between them), the map $\phi_{S^2}^D$ agrees with the map given in Definition \ref{spherical_doubling_map}, justifying the notation. Similarly, when restricted to planar knotoids (knotoidal graphs with three poles and a single segment constituent between two of them), the map $\phi_{S^2}^D$ is precisely the planar reflected doubling map $\phi_{\mathbb{R}^2}^D: \mathbb{K}(\mathbb{R}^2) \to \mathcal{K}(H_3)$ defined in \cite{hypknotoids}.

Note that we can also extend this definition to knotoidal graphs on a general closed orientable projection surface $\Sigma$, by replacing every occurrence of $S^2$ in Definition \ref{defnKGhyp} with $\Sigma$. Again, we obtain a spatial graph $\mathcal{G}$ in a manifold $Y$, and we denote the map associating $(\mathcal{G}, Y)$ to $KG$ by $\phi_{\Sigma}^D$. If $\Sigma$ has genus $g$ and there are no valency-zero poles, then $Y$ is homeomorphic to a thickened closed orientable surface of genus $2g + m - 1$. Otherwise, $Y$ is homeomorphic to a handlebody of genus $(2g + m - 1) + 2(n - m) = 2g + 2n - m - 1$. We can now define hyperbolicity for knotoidal graphs on any closed orientable projection surface $\Sigma$.

\begin{comment}
We remark that one might interpret a knotoidal graph diagram to be a properly immersed spatial graph on a surface with boundary by taking any given knotoidal graph diagram on a surface $\Sigma$ and removing a disk neighborhood of each pole. Then one could eliminate poles entirely from the theory of knotoidal graphs. However, we treat poles as distinct from components of $\partial \Sigma$, so that for instance, the theory of knotoidal graphs on an annulus is markedly distinct from the theory of knotoidal graphs on a sphere. To give a more specific example, we regard a knotoidal graph on an annulus with a single pole as different from a knotoidal graph on a sphere with three poles, even if removing all disk neighborhoods of poles results in knotoidal graphs which are equivalent on thrice-punctured spheres. This distinction is captured by our extended map $\phi_{\Sigma}^D$: a knotoidal graph on an annulus with a single pole is mapped to a spatial graph in a genus 2 handlebody while a knotoidal graph on a sphere with three poles is mapped to a spatial graph in 
\end{comment}

\begin{definition}
Let $KG$ be a knotoidal graph in a closed orientable surface $\Sigma$. Then we say $KG$ is \textit{hyperbolic} if its image $\phi_{\Sigma}^D(KG) = (\mathcal{G}, Y)$ under the $\Sigma$-reflected doubling map is tg-hyperbolic. That is, the manifold $Y \setminus N(\mathcal{G})$ admits a complete hyperbolic metric such that its higher genus boundary components are totally geodesic in the metric. 
\end{definition}

We require the stronger condition of tg-hyperbolicity on $Y \setminus N(\mathcal{G})$ to guarantee that it has a well-defined, finite hyperbolic volume.

\begin{definition}
Let $KG$ be a hyperbolic knotoidal graph in $S^2$, and let $(\mathcal{G}, Y)$ be its image under $\phi_{S^2}^D$ with a tg-hyperbolic metric. Then the \textit{hyperbolic volume} $\vol(KG)$ of $KG$ is defined as
$$
\vol(KG) := \frac{1}{2} \vol(Y \setminus N(\mathcal{G})).
$$
\end{definition}

\begin{figure}[htbp]
    \centering
    \includegraphics[scale = .35]{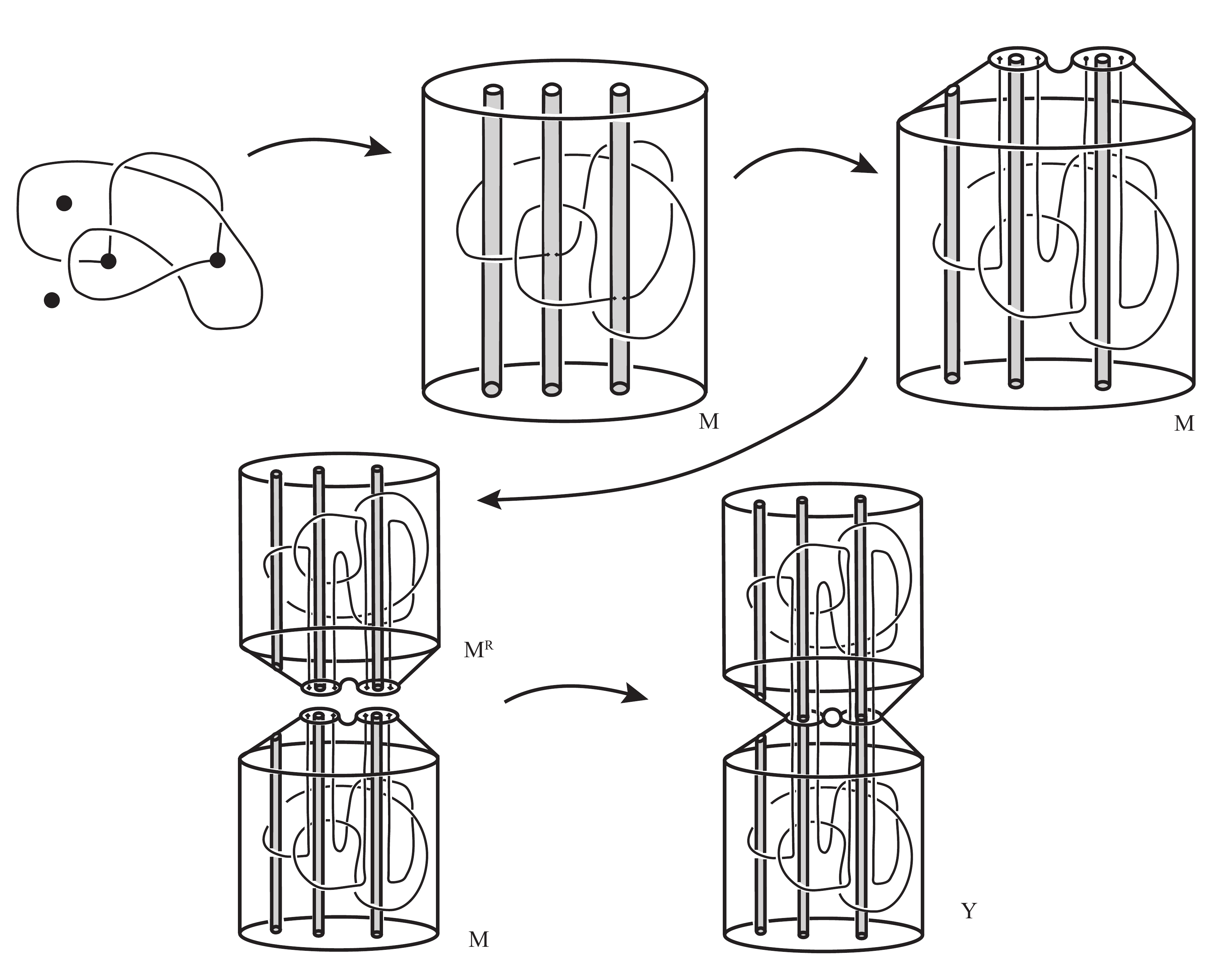}
    \caption{The construction for a generalized knotoid or knotoidal graph.}
    \label{doublingconstruction}
\end{figure}

We remark that hyperbolic volume does not distinguish between a knotoidal graph $KG$ and the knotoidal graph obtained by applying a pole twist move at a pole of $KG$ since the reflected doubling operation ``undoes'' the twist. (However, the polynomial invariants for generalized knotoids defined in Sections \ref{index} and \ref{bracket} may detect twists.)

% {\color{blue} Note that 1-polar knots are never hyperbolic. What happens to 1-staked knots?} {\color{purple} adding a bit about this before the alternating stuff}

\subsection{Staked Links} 
We defined staked links in Example \ref{staked-links}. Given any link $L$ on a closed orientable surface $\Sigma$, we may \textit{stake} it by adding any number of valency-zero poles to its diagram, and we refer to these poles as \textit{stakes}. By means of their rail diagrams, staked links are equivalent to links in handlebodies. 

We would like to consider when staked links are hyperbolic. As a simple example, let $L$ be a link with diagram on $\Sigma = S^2$. Suppose $L'$ is a staked link obtained by adding a single stake to the diagram of $L$. If $(\mathcal{G}, Y)$ is the image of $L'$ under $\phi_{S^2}^D$, then $Y$ is just a 3-ball and $\mathcal{G}$ is a link so that $\mathcal{G}$ in $Y$ is equivalent to $L$ in $S^2 \times I$ under the correspondence between links in $S^3$ and links in $S^2 \times I$. In particular, a 1-staked link $L'$ obtained from $L \subset S^2 \times I$ is hyperbolic if and only if $L$ is hyperbolic. Contrast this with 1-polar links: they are never hyperbolic since there is always an essential annulus in $Y \setminus N(\mathcal{G})$. 

\begin{figure}
    \centering
    \includegraphics[scale = 1.05]{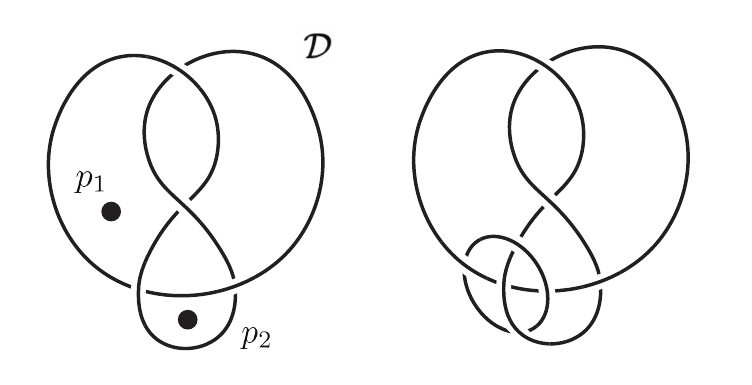}
    \caption{Adding two stakes to a link diagram on $S^2$ is equivalent to adding an unknot component.}
    \label{twostaking}
\end{figure}

As another example, let $L$ be a link with diagram $\mathcal{D}$ on $\Sigma = S^2$, and suppose $L'$ is a staked link obtained by adding two stakes to the diagram of $L$. If $(\mathcal{G}, Y)$ is the image of $L'$ under $\phi_{S^2}^D$, then the complement $Y \setminus N(\mathcal{G})$ is homeomorphic to the complement of a link $L''$ in $S^3$ obtained by adding an unknotted component to $L$ in the following way: let $p_1$ and $p_2$ denote points on $S^2$ where the stakes are added, and add two arcs to the diagram $\mathcal{D}$ with endpoints on $p_1$ and $p_2$ such that one arc crosses over every strand of $\mathcal{D}$ it meets and the other crosses under every strand it meets. This yields a new link diagram $\mathcal{D}''$. Then $L'$ is hyperbolic if and only if $L''$ is hyperbolic. See Figure \ref{twostaking} for an example.

Observe that when $p_1$ and $p_2$ are in distinct, non-adjacent regions, then $L''$ is obtained by augmenting $L$ as in \cite{adams-gen-aug}. When $L$ is a prime, non-split, alternating link which is not a 2-braid, Theorem 2.1 from \cite{adams-gen-aug} shows that such a 2-staked link $L''$ is hyperbolic. 

More generally, if $L$ has a \textit{cellular alternating} diagram on $\Sigma$, then we can determine precisely when a staked link $L'$ obtained by staking the diagram of $L$ is hyperbolic. 

\begin{definition}
Let $\Sigma$ be a closed, orientable surface, and let $L \subset \Sigma \times I$ be a staked link with diagram $\mathcal{D}$. We say $L$ is \textit{weakly prime} if every disk $D \subset \Sigma$ containing at least one crossing of $\mathcal{D}$ such that the circle $\partial D$ intersects $\mathcal{D}$ transversely in two points also contains at least one stake in its interior in a complementary region that does not touch $\partial D$.
\end{definition}

\begin{definition}
Let $\Sigma$ be a closed, orientable surface and let $\mathcal{D}$ be a projection diagram for a link $L$ in $\Sigma$. We say $\mathcal{D}$ is \textit{cellular} if every complementary region $\Sigma \setminus \mathcal{D}$ is a disk. 
\end{definition}

\begin{proposition} \label{alternating-staked-links}
Let $L$ be a link on a closed orientable surface $\Sigma$ and suppose it has a cellular alternating diagram $\mathcal{D}$. Let $L'$ be a staked link obtained by staking $L$, and let $\mathcal{D}'$ be a reduced diagram for $L'$ in $\Sigma$. If $L$ is not the unknot or a $(2, q)$-torus link on $\Sigma = S^2$, then $L'$ is hyperbolic if and only if 
\begin{enumerate}
    \item $\mathcal{D}'$ is weakly prime on $\Sigma$;
    \item every complementary region of $\Sigma \setminus \mathcal{D}'$ contains at most one stake; 
    \item adjacent regions of $\Sigma \setminus \mathcal{D}'$ do not both contain stakes. 
\end{enumerate}
If $L$ is the unknot or a $(2, q)$-torus link on $\Sigma = S^2$, then $L'$ is hyperbolic if and only if conditions (1)-(3) are satisfied and 
\begin{enumerate} \setcounter{enumi}{3}
    \item $L'$ has at least two stakes. 
    \item the rail diagram for $L'$ is not homeomorphic to the complement of a $(p, q)$-torus link in the solid torus, where this  link may be isotoped to sit on the boundary without crossings; 
    \item the rail diagram for $L'$ is not homeomorphic to the thickened torus. 
\end{enumerate}
\end{proposition}

These conditions are all relatively easy to check directly on the diagram $\mathcal{D'}$. Note that if $L$ is the unknot or a $(2, q)$-torus link on $S^2$, then conditions (4) and (5) are automatically satisfied if $L'$ has at least three stakes. Also note that this result subsumes Theorem 2.1 from \cite{adams-gen-aug}. 

\begin{proof}
As in the remarks following the definitions of generalized knotoids and knotoidal graphs, adding a stake to a link diagram $\mathcal{D}$ on $\Sigma$ is equivalent to removing an open disk neighborhood of $\Sigma$. Hence, we can view $\mathcal{D}'$ as a link diagram on a compact surface $\Sigma'$ with boundary. Then the proposition is a restatement of Theorem 1.6 and Corollary 1.7 %{\color{purple} warning: make sure to change number} 
from \cite{alternatingboundary}. 
\end{proof}

This proposition gives us a large class of hyperbolic staked links. The following theorem uses it to show that an even larger class is hyperbolic.

\begin{theorem} \label{checkerboard-staking}
Let $\Sigma$ be a closed, orientable surface and let $L \subset \Sigma \times I$ be a link with a cellular projection diagram $\mathcal{D}$ with at least two crossings such that $\Sigma \setminus \mathcal{D}$ is checkerboard-colorable. If $\mathcal{D}$ is not a diagram of the 2-component unlink or Hopf link on $S^2$ with exactly two crossings, then $\mathcal{D}$ can be staked such that the resulting staked link is hyperbolic. %{\color{purple} added cellular, also exceptional cases}
\end{theorem}

\begin{proof}
Consider a checkerboard coloring of $\Sigma \setminus \mathcal{D}$ and in every shaded face, add a single stake. We call this a \textit{checkerboard staking} of $\mathcal{D}$, and we claim that the resulting staked link is hyperbolic. We may view it as a link in the thickened surface $\Sigma'$ with boundary obtained by removing an open disk neighborhood of each stake in $\Sigma$. By abuse of notation, we refer to the link viewed in the handlebody $\Sigma' \times I$ by $L$. We want to show that $M = (\Sigma' \times I) \setminus N(L)$ is tg-hyperbolic. 

\begin{figure} [htbp]
    \centering
    \includegraphics[scale = 0.85]{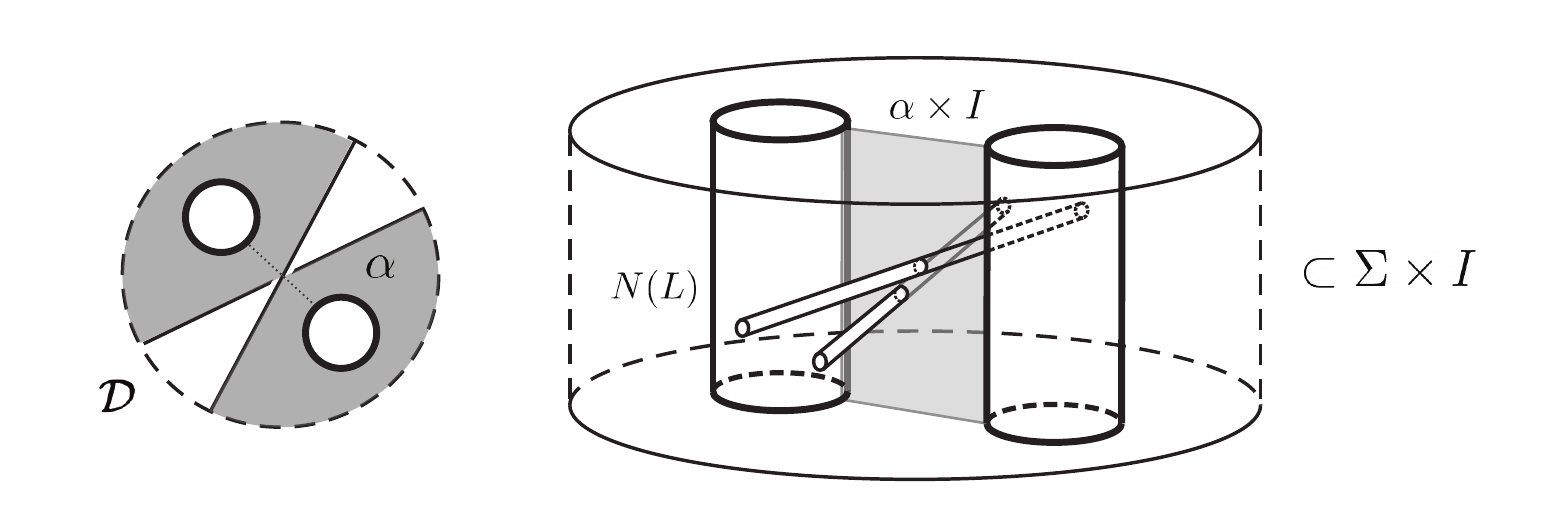}
    \caption{A checkerboard staking near a crossing; crossing the arc $\alpha$ by $I$ yields the twice-punctured disk $F = \alpha \times I$ in $\Sigma \times I$, shown on the right.} 
    \label{localcrossingpic}
\end{figure}

Each crossing of the resulting diagram $\mathcal{D}'$ locally divides $\Sigma'$ into four regions, and two of the regions opposite to each other contain a circle boundary where a neighborhood of the stake was removed. For some crossing $c$, denote these regions by $R_1$ and $R_2$, and denote the respective circle boundaries by $C_1$ and $C_2$. Let $\alpha$ be an arc in $\Sigma'$ with endpoints on the $C_i$ such that it intersects $\mathcal{D}$ exactly once through the crossing $c$. Then $\alpha \times I$ is a properly embedded disk $F$ in $\Sigma' \times I$ which is twice-punctured by $L'$ at the crossing. See Figure \ref{localcrossingpic}. 

 We can cut along $F$, yielding copies $F_1$ and $F_2$, rotate copy $F_2$ by $2 \pi$, and reglue the two copies back together, changing the crossing to the opposite crossing. This is a homeomorphism of $M$ so the resulting link in the handlebody $\Sigma' \times I$ is hyperbolic if and only if the original link is hyperbolic. We may switch any of the crossings of $L$ in the handlebody while preserving hyperbolicity. Therefore, we can perform this move several times to obtain an alternating link $L'$ in a handlebody. 

Now, we apply Proposition \ref{alternating-staked-links}.
Conditions (2) and (3) are automatically satisfied by construction. If $\Sigma = S^2$, then we need to consider conditions (4)-(6). Condition (4) is satisfied since $\mathcal{D}$ has at least two crossings: hence, there are at least four regions of $\Sigma \setminus \mathcal{D}$. Note that conditions (5) and (6) are automatically satisfied if there are at least three stakes. We can choose such a checkerboard staking if there are at least five regions; hence, we may assume that there are exactly two crossings and four regions of $\Sigma \setminus \mathcal{D}$ such that exactly two regions are shaded and two are white in a checkerboard coloring. This only occurs if $\mathcal{D}$ is the diagram of a 2-component unlink or Hopf link wth exactly two crossings. 

It remains to check condition (1). We may assume that $\mathcal{D}$ is a reduced diagram on $\Sigma'$. By definition, it suffices to check that every circle $\gamma \subset \Sigma$ which bounds a disk $D$ in $\Sigma$ containing crossings of $\mathcal{D}$ and is twice-punctured by $\mathcal{D}$ contains a region of each color. Since $D$ contains at least one crossing, it must contain at least one region. Without loss of generality, say this region is shaded. If it is staked, then we are done; otherwise, $\mathcal{D}$ is staked in the unshaded regions. Since $\gamma$ meets $\mathcal{D}$ exactly twice, it intersects the shaded faces in a single arc and the unshaded faces in another arc. Let $W$ be the unshaded face intersecting $\gamma$; since $\gamma$ contains only shaded faces, $W$ is the only unshaded face $\gamma$ runs through. Let $c$ be a crossing contained in $D$ and let $\mu$ be a small arc in $W \cup \mathcal{D}$ which intersects $\mathcal{D}$ exactly once through the double point of $c$. Since $W$ is connected, there is an arc $\sigma \subset W$ connecting the endpoints of $\mu$. Then $\mu \cup \gamma$ is a circle which bounds a disk in $\Sigma'$ and intersects $\mathcal{D}$ precisely once through a double point. This contradicts $\mathcal{D}'$ being reduced on $\Sigma'$. Hence, all of the conditions of Proposition \ref{alternating-staked-links} are satisfied and the theorem follows. 

\begin{comment}
We can ensure that $L'$ has at least two stakes by adding arbitrarily many Type I Reidemeister moves to the original link diagram $\mathcal{D}$ before checkerboard coloring {\color{purple} we want to say original diagram can be staked}. To check Condition (1), by the definition of weakly prime it suffices to check that every circle $\gamma \subset F$ containing crossings of $\mathcal{D}$ contains a face of each color. Then, we are guaranteed to have a pole within $\gamma$, so that $\ell$ does not bound a disk in $F$. 
\end{comment}
\end{proof}

\begin{corollary} \label{every-link-sphere-stake}
Every link $L$ in $S^2 \times I$ has a diagram which can be staked such that the resulting staked link is hyperbolic. 
\end{corollary}

\begin{proof}
Every link diagram in $S^2$ is checkerboard-colorable. Also, every link diagram in $S^2$ can be made connected with sufficiently many crossings (via Type I and II Reidemeister moves) so that the previous theorem applies.
\end{proof}

The checkerboard staking is maximal in the sense that adding an additional stake to $\mathcal{D}$ anywhere in the diagram causes it to no longer be hyperbolic. There are examples of links in $\Sigma$ which are not hyperbolic prior to staking but become hyperbolic after adding just a few stakes. For instance, if $L$ becomes alternating after switching crossings $c_1, \dotsc, c_n$, then we may repeat the proof of Theorem \ref{checkerboard-staking} after staking once in each shaded region adjacent to one of the crossings $c_i$. In particular, if $L$ is a link which is \textit{almost alternating}, that is, becomes alternating after switching just one crossing, then it can be made hyperbolic by staking just twice (assuming the other conditions in the statement of Proposition \ref{alternating-staked-links} are satisfied). 

\begin{figure}
    \centering
    \includegraphics[scale=0.9]{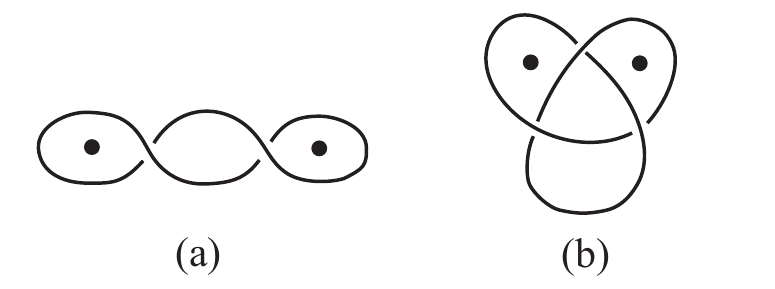}
    \caption{Staked diagrams of the unknot and trefoil realizing the minimum volume over all staked diagrams.} 
    \label{staked-unknot-trefoil}
\end{figure}

However, not all hyperbolic staking configurations arise this way. For example, see the  the staking configurations of the unknot and trefoil in Figure \ref{staked-unknot-trefoil}.  
%{\color{red} The staking from the infinite family is three stakes, and I don't think its the canonical staking you describe here, do we still want to include it?} {\color{purple} I guess here a reason to keep it is to show it's not just a coincidence of small numbers? like the unknot and trefoil might have weird staking configurations because they happen to be problematic knots in view of the thm statement. But I'm also fine taking it out.}  
Hence, given a link $L$, it would be interesting to determine exactly the  staking configurations that yield a hyperbolic staked link, that is, extend Proposition \ref{alternating-staked-links} to larger classes of links.

\begin{definition}
We define the \emph{staked volume} of a link $L$ in $\Sigma \times I$, where $\Sigma$ is a closed surface, to be the minimum volume over all diagrams $\mathcal{D}$ of $L$ in $\Sigma \times \{1/2\}$ of staked links obtained by adding stakes to $\mathcal{D}$. If no such staking exists for any diagram $\mathcal{D}$, then we set the staked volume of $L$ to infinity. 
\end{definition}

By Corollary \ref{every-link-sphere-stake}, we may associate a finite hyperbolic volume to every link in $S^3$, including links that are not hyperbolic in the classical sense. In this regard, staked volume is similar to Turaev volume as defined in \cite{turaevvolume}. 

\begin{example} \label{stakedvolunknot}
The staked volume of the trivial knot in $S^3$ is  $ 3.6638 \dots$, and it is achieved by the staking configuration in Figure \ref{staked-unknot-trefoil}(a). To see this, observe that for any staked link on $S^2$ with exactly two stakes, its volume as a knotoidal graph is computed as the volume of a link complement in the solid torus. In particular, this is a 2-cusped manifold. In \cite{min2cusped}, Agol proved that the minimum volume of a 2-cusped manifold is $ 3.6638 \dots$. Hence, this must be the minimum volume of any staked diagram of the trivial knot with two stakes. Any diagram of the trivial knot with one stake is isotopic to the trivial diagram with one stake, which yields a trivial knot in the 3-sphere, which is not hyperbolic.  Any diagram of the trivial knot with three or more stakes yields a knot in a handlebody of genus at least two. However, for this to be tg-hyperbolic, the boundary must be totally geodesic and by Theorem 5.2 of \cite{miyamoto}, we know that any manifold with totally geodesic boundary must have volume at least $3.6638\dots$ and if it has volume exactly $3.6638\dots$, it cannot be the complement of a knot or link in a handlebody. 
\end{example}

% Any staking configuration of the unknot with three or more stakes yields a knot complement in a handlebody with genus at least 2. To compute the volume, we double across the higher genus boundary, yielding a 2-cusped manifold. {\color{purple} run into an issue with the volume bound...possible that when volume is halved, it dips below 3.6638??}

\begin{example}
The staked volume of the trefoil knot in $S^3$ is also $3.6638 \dots$, and it is achieved by the staking configuration in Figure \ref{staked-unknot-trefoil}(b). Similarly to Example \ref{stakedvolunknot}, it is the minimum volume of any staked diagram of the trefoil with exactly two stakes.  As in the previous example, one stake will not make it hyperbolic and three or more stakes will cause it volume to be larger. 
\end{example}

\begin{example}
For any twist knot other than the trefoil knot, its staked volume is its hyperbolic volume. If we use one stake, then removing the corresponding rail from $S^2 \times I$ yields the link in a ball, and capping that ball with another ball yields the knot in $S^3$, with volume the original hyperbolic volume of the knot. If we utilize two stakes, the result will be the complement of a 2-component link in $S^3$, which by \cite{min2cusped} has volume at least $3.6638\dots$. In fact there is a 2-pole staking for each twist knot that yields this volume. However, every twist knot comes from Dehn filling the Whitehead link, which has this volume, and Dehn filling always lowers volume, so the volume from one stake is strictly less than the volume from two stakes. And exactly as in the last two examples, the volume from three or more  stakes is strictly greater than $3.6638\dots$.
\end{example}

%commenting this out since our new theorem about staked knots is much more powerful; we can always add it back in though!
% Am adding back in for the almost alternating staked example
%\begin{figure}[htpb]
    %\centering
  %  \includegraphics[scale=0.35]{images/staked_example_after_stakes.pdf}
%\caption{An example of an almost alternating family of knots which are not hyperbolic without stakes as they are composite knots, but are hyperbolic after adding three stakes as shown. Here $T,T'$ are prime, alternating tangles which have closures that are not $2$-braids; it is proved in \cite{hypknotoids} that this staked family is hyperbolic.}
%\end{figure}

\bibliographystyle{plain}
\bibliography{references}

\begin{thebibliography}{10}

\bibitem{adams-gen-aug}
C.~Adams.
\newblock Generalized augmented alternating links and hyperbolic volumes.
\newblock {\em Algebraic \& Geometric Topology}, 17(6):3375--3397, Oct 2017.

\bibitem{hypknotoids}
C.~Adams, A.~Bonat, M.~Chande, J.~Chen, M.~Jiang, Z.~Romrell, D.~Santiago,
  B.~Shapiro, and D.~Woodruff.
\newblock Hyperbolic knotoids.
\newblock {\em preprint}, 2022.

\bibitem{alternatingboundary}
C.~Adams and J.~Chen.
\newblock Hyperbolicity of alternating links in thickened surfaces with
  boundary.
\newblock {\em in preparation}, 2022.

\bibitem{turaevvolume}
C.~Adams, O.~Eisenberg, K.~Kapoor J.~Greenberg, Z.~Liang, K.~O'Connor,
  N.~Pacheco-Tallaj, and Y.~Wang.
\newblock Turaev hyperbolicity of classical and virtual knots.
\newblock {\em Algebraic and Geometric Topology}, 21:3459--3482, 2021.

\bibitem{min2cusped}
I.~Agol.
\newblock The minimal volume orientable hyperbolic 2-cusped 3-manifolds.
\newblock {\em Proceedings of the American Mathematical Society},
  138(10):3723--3732, May 2010.

\bibitem{winding-poly}
K.~Bhataineh.
\newblock An invariant of planar knotoids and finite-type invariants, 2020.
\newblock Preprint available at
  \url{https://www.researchgate.net/publication/344446254_AN_INVARIANT_OF_PLANAR_KNOTOIDS_AND_FINITE-TYPE_INVARIANTS}.

\bibitem{polar}
K.~Bhataineh.
\newblock New polynomial invariants of knotoids and the theory of polar knots.
\newblock {\em Mediterranean Journal of Mathematics}, 19:Article 40, 2022.

\bibitem{boden-rushworth}
H.~U. Boden and W.~Rushworth.
\newblock Minimal crossing number implies minimal supporting genus.
\newblock {\em Bulletin of the London Mathematical Society}, 53:1174--1184,
  2021.

\bibitem{knotoidproteins}
J.~Dorier, D.~Goundaroulis, F.~Benedetti, and A.~Stasiak.
\newblock Knotoid: a tool to study the entanglement of open protein chains
  using the concept of knotoids.
\newblock {\em Bioinformatics}, 34(19):3402--3404, 2018.

\bibitem{wriggle}
L.~Folwaczny and L.~Kauffman.
\newblock A linking number definition of the affine index polynomial and
  applications.
\newblock {\em Journal of Knot Theory and its Ramifications}, 22:1341004 (30
  pages), 2013.

\bibitem{multi-linkoid}
B.~Gabrovšek and N.~Gügümcü.
\newblock Invariants of multi-linkoids.
\newblock {\em arXiv:2204.11234}, 2022.

\bibitem{bondedknotoids}
D.~Goundaroulis, N.~Gügümcü ande S.~Lambropoulou, J.~Dorier, A.~Stasiak, and
  L.~Kauffman.
\newblock Topological models for open-knotted protein chains using the concepts
  of knotoids and bonded knotoids.
\newblock {\em Polymers}, 9:444, 2017.

\bibitem{knotoidproteins2}
D.~Goundaroulis, J.~Dorier, and A.~Stasiak.
\newblock Knotoids and protein structure.
\newblock {\em Topol. Geom. Biopolym.}, 746, 2020.

\bibitem{bondedknotoids2}
N.~Gügümcü, B.~Gabrovsek, and L.H. Kauffman.
\newblock Invariants of bonded knotoids and applications to protein folding.
\newblock {\em Symmetry}, 14, 2022.

\bibitem{NIV}
N.~Gügümcü and L.~Kauffman.
\newblock New invariants of knotoids.
\newblock {\em European Journal of Combinatorics}, 65:186--229, 2017.

\bibitem{gp-min-crossing}
N.~Gügümcü and L.~Kauffman.
\newblock Parity, virtual closure and minimality of knotoids.
\newblock {\em Journal of Knot Theory and its Ramifications}, 30:2150076 (28
  pages), 2022.

\bibitem{graphoids}
N.~Gügümcü, L.~Kauffman, and P.~Pongtanapaisan.
\newblock Graphoids.
\newblock {\em preprint to appear on the ArXiv}, 2022.

\bibitem{henrich}
A.~Henrich.
\newblock A sequence of degree-one vassiliev invariants for virtual knots.
\newblock {\em Journal of Knot Theory and its Ramifications}, 19:461--487,
  2010.

\bibitem{Kashaev}
R.~Kashaev.
\newblock Invariants of long knots.
\newblock In {\em Representation Theory, Mathematical Physics, and Integrable
  Systems}, pages 431--451. Birkh\"auser, 2019.

\bibitem{affine-index-virtual}
L.~Kauffman.
\newblock An affine index polynomial invariant of virtual knots.
\newblock {\em Journal of Knot Theory and its Ramifications}, 22:1340007 (30
  pages), 2013.

\bibitem{bonded}
L.~Kauffman, B.~Gabrovšek, and N.~Gügümcü.
\newblock Topological invariants of bonded proteins, 2022.
\newblock Preprint available at
  \url{https://www.researchgate.net/publication/361633189_TOPOLOGICAL_INVARIANTS_OF_BONDED_PROTEINS}.

\bibitem{family-index-poly-knotoid}
S.~Kim, Y.~H. Im, and S.~Lee.
\newblock A family of polynomial invariants for knotoids.
\newblock {\em Journal of Knot Theory and its Ramifications}, 27:1843001 (15
  pages), 2018.

\bibitem{rails}
D.~Kodokostas and S.~Lambropoulou.
\newblock Rail knotoids.
\newblock {\em Journal of Knot Theory and its Ramifications}, 28:1940019 (19
  pages), 2019.

\bibitem{kutluay}
D.~Kutluay.
\newblock Winding homology of knotoids, 2020.
\newblock Ph.D. thesis; arXiv:2002.07871.

\bibitem{linov}
L.~Linov.
\newblock Signed heights of knotoids.
\newblock {\em Journal of Knot Theory and its Ramifications}, 31:2250037 (25
  pages), 2022.

\bibitem{manturovparity}
V.~O. Manturov.
\newblock Parity and projection from virtual knots to classical knots.
\newblock {\em Journal of Knot Theory and its Ramifications}, 22:1350044, 2013.

\bibitem{miyamoto}
Yosuke Miyamoto.
\newblock Volumes of hyperbolic manifolds with totally geodesic boundary.
\newblock {\em Topology}, 33(4):613--629, 1994.

\bibitem{HL}
X.-S.~Lin N.~Habegger.
\newblock The classification of links up to link-homotopy.
\newblock {\em J. Amer. Math. Soc.}, 3:389–419, 1990.

\bibitem{Turaev2}
V.~Turaev.
\newblock Cobordism of knots on surfaces.
\newblock {\em Journal of Topology}, 1:285--305, 2008.

\bibitem{Turaev1}
V.~Turaev.
\newblock Knotoids.
\newblock {\em Osaka Journal of Mathematics}, 49:195--223, 02 2010.

\end{thebibliography}

\end{document}